\documentclass[11pt]{article}
\usepackage{amsfonts,mathrsfs,amssymb,amsthm,mathptm}
\usepackage{amsmath,amscd}
\usepackage{mathptm,pslatex}
\usepackage{multirow}
\usepackage{color}
\usepackage[all]{xy}
\usepackage{url}
\usepackage{cite}
\usepackage{exscale}
\usepackage{relsize}
\usepackage{bbm}
\begin{document}
\newtheorem{Def}{Definition}[section]
\newtheorem{Bsp}[Def]{Example}
\newtheorem{Prop}[Def]{Proposition}
\newtheorem{Theo}[Def]{Theorem}
\newtheorem{Lem}[Def]{Lemma}
\newtheorem{Koro}[Def]{Corollary}
\theoremstyle{definition}
\newtheorem{Rem}[Def]{Remark}

\newcommand{\add}{{\rm add}}
\newcommand{\Coker}{{\rm Coker}}
\newcommand{\Cone }{{\rm Cone}}
\newcommand{\gd}{{\rm gl.dim}}
\newcommand{\dm}{{\rm dom.dim}}
\newcommand{\DTr}{{\rm DTr}}
\newcommand{\E}{{\rm E}}
\newcommand{\Mor}{{\rm Morph}}
\newcommand{\End}{{\rm End}}
\newcommand{\Ext}{{\rm Ext}}
\newcommand{\Fac}{{\rm Fac}}
\newcommand{\fd} {{\rm fin.dim}}
\newcommand{\fld}{{\rm flat.dim}}
\newcommand{\Gen}{{\rm Gen}}
\newcommand{\Hom}{{\rm Hom}}
\newcommand{\holim}{{\rm Holim}}
\newcommand{\hocolim}{{\rm Hocolim}}
\newcommand{\id}{{\rm inj.dim}}
\newcommand{\ind}{{\rm ind}}
\newcommand{\Img  }{{\rm Im}}
\newcommand{\Ker  }{{\rm Ker}}
\newcommand{\ol}{\overline}
\newcommand{\overpr}{$\hfill\square$}
\newcommand{\pd}{{\rm proj.dim}}
\newcommand{\projdim}{\pd}
\newcommand{\rad}{{\rm rad}}
\newcommand{\rsd}{{\rm res.dim}}
\newcommand{\rd} {{\rm rep.dim}}
\newcommand{\soc}{{\rm soc}}
\renewcommand{\top}{{\rm top}}

\newcommand{\cpx}[1]{#1^{\bullet}}
\newcommand{\D}[1]{{\mathscr D}(#1)}
\newcommand{\Dz}[1]{{\mathscr D}^+(#1)}
\newcommand{\Df}[1]{{\mathscr D}^-(#1)}
\newcommand{\Db}[1]{{\mathscr D}^b(#1)}
\newcommand{\C}[1]{{\mathscr C}(#1)}
\newcommand{\Cz}[1]{{\mathscr C}^+(#1)}
\newcommand{\Cf}[1]{{\mathscr C}^-(#1)}
\newcommand{\Cb}[1]{{\mathscr C}^b(#1)}
\newcommand{\K}[1]{{\mathscr K}(#1)}
\newcommand{\Kz}[1]{{\mathscr K}^+(#1)}
\newcommand{\Kf}[1]{{\mathscr  K}^-(#1)}
\newcommand{\Kb}[1]{{\mathscr K}^b(#1)}
\newcommand{\modcat}{\ensuremath{\mbox{{\rm -mod}}}}
\newcommand{\Modcat}{\ensuremath{\mbox{{\rm -Mod}}}}
\newcommand{\gm}{{\rm _{\Gamma_M}}}
\newcommand{\gmr}{{\rm _{\Gamma_M^R}}}
\newcommand{\stmodcat}[1]{#1\mbox{{\rm -{\underline{mod}}}}}
\newcommand{\StHom}{{\rm \underline{Hom}}}
\newcommand{\pmodcat}[1]{#1\mbox{{\rm -proj}}}
\newcommand{\imodcat}[1]{#1\mbox{{\rm -inj}}}
\newcommand{\Pmodcat}[1]{#1\mbox{{\rm -Proj}}}
\newcommand{\Imodcat}[1]{#1\mbox{{\rm -Inj}}}
\newcommand{\opp}{^{\rm op}}
\newcommand{\otimesL}{\otimes^{\rm\mathbb L}}
\newcommand{\rHom}{{\rm\mathbb R}{\rm Hom}\,}

\def\vez{\varepsilon}\def\bz{\bigoplus}  \def\sz {\oplus}
\def\epa{\xrightarrow}
\def\inja{\hookrightarrow}

\newcommand{\lra}{\longrightarrow}
\newcommand{\lraf}[1]{\stackrel{#1}{\lra}}
\newcommand{\ra}{\rightarrow}
\newcommand{\dk}{{\rm dim_{_{k}}}}
\newcommand{\colim}{{\rm colim\, }}
\newcommand{\limt}{{\rm lim\, }}
\newcommand{\Add}{{\rm Add }}
\newcommand{\Prod}{{\rm Prod }}
\newcommand{\Tor}{{\rm Tor}}
\newcommand{\Cogen}{{\rm Cogen}}

{\Large \bf
\begin{center}
Symmetric subcategories, tilting modules and derived recollements
\end{center}}

\bigskip
\centerline{\textbf{Hongxing Chen} and \textbf{Changchang Xi}$^*$}

\renewcommand{\thefootnote}{\alph{footnote}}
\setcounter{footnote}{-1} \footnote{ $^*$ Corresponding author.
Email: xicc@cnu.edu.cn; Fax: 0086 10 68903637.}
\renewcommand{\thefootnote}{\alph{footnote}}
\setcounter{footnote}{-1} 
\renewcommand{\thefootnote}{\alph{footnote}}
\setcounter{footnote}{-1} 

\begin{abstract}
For any good tilting module $T$ over a ring $A$, there exists an $n$-symmetric subcategory $\mathscr{E}$ of a module category such that the derived category of the endomorphism ring of $T$ is a recollement of the derived categories of $\mathscr{E}$ and $A$ in the sense of Beilinson-Bernstein-Deligne. Thus the kernel of the total left-derived tensor functor induced by the tilting module is triangle equivalent to the derived category of $\mathscr{E}$.
\end{abstract}

{\footnotesize\tableofcontents}

\section{Introduction \label{sect1}}
In representation theory of algebras and groups, finitely generated tilting modules have successfully been applied to understanding different aspects of algebraic and homological structure and properties of finite groups, finite-dimensional algebras and modules (for example, see \cite{BBu, CPS, happel87, hr, Ri}), while infinitely generated tilting modules not only have involved many important modules, such as adic modules, Fuchs divisible modules, generic modules, and Pr\"ufer modules, but also have been of significant interest in studying derived categories and equivalences of general algebras and rings (see \cite{AC, Bz1, Bz2, XC1, XC2,yang}). In this general context, a lot of results emerge with completely different features. So, the theory of infinitely generated tilting modules renews our view on finitely generated tilting modules, and gives us also surprising information about the whole tilting theory.

Let us first recall the definition of tilting modules.

\begin{Def}{\rm \cite{BBu, hr, miy, ct, AC}} \label{Tilting} Let $n\ge 0$ be a natural number and $A$ a unitary ring. A
left $A$-module $T$ is called an \emph{$n$-tilting} $A$-module
if the three conditions hold:

$(T1)$ There is an exact sequence of $A$-modules:
$$ 0\lra P_n\lra \cdots\lra P_1\lra P_0\lra {}_AT\lra 0$$
with all $P_i$ projective $A$-modules, that is, $\pd(_AT)\le n$.

$(T2)$ For all $j\geq 1$ and nonempty sets
$\alpha$, $\Ext^j_A(T, T^{(\alpha)})=0$, where $T^{(\alpha)}$ stands for the direct sum of $\alpha$ copies of $T$.

$(T3)$ There exists an exact sequence of $A$-modules
$$ 0\lra {}_AA \lra T_0 \lra
T_1\lra \cdots \lra T_n\lra 0,$$ where $T_i$ is
isomorphic to a direct summand of a direct sum of copies of
$T$ for $ 0\leq i\leq n$.
\end{Def}

Recall that an $n$-tilting $A$-module $T$ is said to be \emph{good} (see \cite{Bz2}) if each $T_i$ in $(T3)$ is
isomorphic to a direct summand of the direct sum of finitely many copies of
$T$.

Given an arbitrary tilting module $T'$, one can always find a good tilting module $T$  such that $T$ and $T'$ generate the same full subcategories in $A\Modcat $, the category of all left $A$-modules, that is, $T$ and $T'$ are equivalent (see \cite{Bz2}).

For an exact category $\mathscr{E}$, we denote by $\D{\mathscr{E}}$ the unbounded derived category of $\mathscr{E},$ and for a ring $A$, we denote by $\D{A}$ the unbounded derived module category of $A$, that is, $\D{A}=\D{A\Modcat}.$

The purpose of this article is to prove the following derived recollement theorem for arbitrary good tilting modules over rings. For the definition of $n$-symmetric subcategories, we refer to Section \ref{sect3}.

\begin{Theo}\label{Main result} Let $A$ be a ring, $T$ a
good $n$-tilting $A$-module, and $B$ the endomorphism ring of $_AT$. Then there exist an $n$-symmetric subcategory $\mathscr{E}$ of $B\Modcat$ and  a recollement of unbounded derived categories

$$\xymatrix{\D{\mathscr{E}}\ar[r]&\D{B}\ar[r]^{G}\ar@/^1.2pc/[l]\ar@/_1.2pc/[l]
&\D{A}\ar@/^1.2pc/[l]\ar@/_1.2pc/[l]}\vspace{0.2cm}$$
where $G$ is the total left-derived tensor functor defined by $_AT_B$. Moreover, this recollement restricts to a recollement of bounded-above derived categories.
\end{Theo}

Actually, for a tilting $A$-module $T$ of projective dimension $n$, the $n$-symmetric subcategory in Theorem \ref{Main result} can precisely be described as
$$\mathscr{E}:=\{Y\in B\Modcat \mid \Tor_i^B(T,Y)=0\mbox{\;for\, all \,} i\geq0\}.$$

Recall that Happel's Theorem (\cite{happel87}, see also \cite{CPS}) for finitely generated tilting modules ensures a triangle equivalence of the derived categories of the rings $A$ and $B$, while for an arbitrary infinitely generated tilting module, one cannot get a derived equivalence but a recollement of triangulated categories in which one category seems difficult to understand (see \cite{Bz2} or Theorem \ref{Torsion}). In \cite{XC1, XC2} special efforts are made to describe this category as a derived module category for good tilting modules with special properties. As pointed in \cite{XC2}, this does not always happen. So it remains mysterious what a general form of Happel's Theorem looks like for arbitrary infinitely generated tilting modules. Now, Theorem \ref{Main result} may serve as a counterpart of Happel's Theorem for infinitely generated good tilting modules. Furthermore, compared with the main results in \cite{yang} and \cite{BaPa}, Theorem \ref{Main result} looks more intrinsic and the derived category of the symmetric subcategory $\mathscr{E}$ seems much simpler than the derived category of a differential graded algebra.

A crucial point in the proof of Theorem \ref{Main result} is to show that the kernel of the total left-derived tensor functor ${_A}T\otimesL_B-$ can be realised as the derived category of the $n$-symmetric subcategory $\mathscr{E}$. When $\mathscr{E}$ is an abelian subcategory,
this can be shown by analysing cohomologies of a special complex of $B$-modules which
is a compact generator for the kernel (for example, see \cite{HKL,XC1,BaPa,XC2}). In general,
however, $\mathscr{E}$ is not always an abelian subcategory. Thus the methods developed in the foregoing-mentioned papers do not work anymore because cohomologies of complexes over $\mathscr{E}$
do not necessarily belong to $\mathscr{E}$.
To circumvent this obstacle here, we introduce two triangle functors from $\D{B}$ to $\D{\mathscr{E}}$ and from $\D{B}$ to $\D{\mathscr{E}^{\bot}}$ (see Section \ref{sect4.1} for more details), where $\mathscr{E}^{\bot}$ is the right perpendicular subcategory to $\mathscr{E}$ in $B\Modcat$. We then realise the kernel of  ${_A}T\otimesL_B-$ as derived category of the symmetric subcategory $\mathscr{E}$ of $B$-Mod through a bridge of a $t$-structure induced from the tilting module.
Notably, the method in this article sheds  some new light on when the kernel of  ${_A}T\otimesL_B-$ is triangle equivalent to the derived module category of a ring, that is, when a good titling module $T$ is homological (see \cite{XC2}). Also, the method establishes a connection of homological tilting modules with derived decompositions of module categories (see \cite{XC5} or Section \ref{sect3.3} for Definition). We may state these observations as a corollary.

\begin{Koro}\label{Homological}
The following are equivalent for a good tilting $A$-module $T$:

$(1)$ $T$ is a homological tilting module.

$(2)$ $\mathscr{E}$ is an abelian subcategory of $B\Modcat$.

$(3)$ $H^m(\Hom_A(\cpx{P}, A)\otimes_BT)=0$ for all $m\geq 2$, where
$\cpx{P}$ is a deleted projective resolution of ${_A}T$, and $H^m$ is the $m$-th cohomology functor.

$(4)$ $(\mathscr{E}, \mathscr{E}^{\bot})$ is a derived decomposition of the abelian category $B\Modcat$, where $\mathscr{E}^{\bot}:=\{Y\in B\Modcat\mid \Ext^n_B(X,Y)=0, \forall \; X\in \mathscr{E}, n\ge 0\}$.
\end{Koro}
The equivalence of $(1)$ and $(2)$ is known in \cite[Proposition 6.2]{BaPa}, while the equivalence of (1) and (3)
is proved in \cite[Theorem 1.1]{XC2}. Only the condition (4) is new.

Finally, we point out that Theorem \ref{Main result} can be applied to special rings to get recollemnts of bounded derived categories.

\begin{Koro} If $A$ is a left coherent ring (that is, every finitely generated left ideal of $A$ is finitely presented) and $T$ is a good tilting $A$-module with $B:=\End_A(T)$, then, for $*\in \{b,+,-,\emptyset \}$,  there exists a recollement of derived categories
$$\xymatrix{\mathscr{D}^*(\mathscr{E})\ar[r]&\mathscr{D}^*(B)\ar[r]^{G}\ar@/^1.2pc/[l]\ar@/_1.2pc/[l]
&\mathscr{D}^*(A)\ar@/^1.2pc/[l]\ar@/_1.2pc/[l]}\vspace{0.2cm}$$
where $\mathscr{E}$ is a symmetric subcategory of $B\Modcat$ and $G$ is the total left-derived tensor functor defined by $_AT_B$.
\label{Koro1.4}
\end{Koro}

This article is organised as follows. In Section \ref{sect2}, we recall definitions, notations and
basic facts needed for proofs. In Section \ref{sect3}, we introduce $m$-symmetric subcategories for each natural number $m$ and give methods to construct such subcategories. In Section \ref{sect4}, we prove the main result, Theorem \ref{Main result},  and its corollaries.

\medskip
{\bf Acknowledgement.}
The research work of both authors was partially supported by the National Natural Science Foundation of China (12031014) and the Beijing Natural Science Foundation (1192004).

\section{Preliminaries\label{sect2}}
In this section, we briefly recall some definitions, basic facts and
notation used in this paper. For unexplained notation employed in
this paper, we refer the reader to \cite{XC1,XC2} and the references
therein.

\subsection{Semi-orthogonal decompositions, recollements and homotopy (co)limits \label{sect2.1}}

Let $\mathcal C$ be an additive category.

A full subcategory of $\mathcal
C$ is always assumed to be closed under isomorphisms.
For an object $X\in \mathcal{C}$, we write $\add(X)$ for the full
subcategory of $\mathcal{C}$ consisting of all direct summands of
finite coproducts of copies of $X$. If $\mathcal{C}$ admits
coproducts (that is, coproducts indexed over sets exist in
${\mathcal C}$), we write $\Add(X)$ for the full subcategory of
$\mathcal{C}$ consisting of all direct summands of coproducts
of copies of $X$. Dually, we write $\Prod(X)$ for the full subcategory of $\mathcal{C}$
consisting of all direct summands of products of copies of
$X$ if $\mathcal{C}$ admits products (that is, products indexed over sets exist in
${\mathcal C}$).

We denote the composition of two morphisms $f: X\to Y$ and $g: Y\to Z$ in $\mathcal C$ by $fg$. We write $f^*$ for $\Hom_{\mathcal
C}(Z,f):\Hom_{\mathcal C}(Z,X)\ra \Hom_{\mathcal C}(Z,Y)$ and $f_*$
for $\Hom_{\mathcal C}(f,Z): \Hom_{\mathcal C}(Y, Z)\ra \Hom_{\mathcal
C}(X, Z)$. While for two functors $F:\mathcal {C}\to
\mathcal{D}$ and $G: \mathcal{D}\to \mathcal{E}$, the composition of $F$ and $G$ is denoted by $GF$ which is a functor from $\mathcal C$ to
$\mathcal E$. The kernel and image of $F$ are defined as $\Ker(F):=\{X\in \mathcal{C}\mid FX\simeq 0\}$ and
$\Img(F):=\{Y\in \mathcal{D}\mid \exists\, X\in \mathcal{C},\; FX\simeq Y\}$, respectively.
Thus $\Ker(F)$ and $\Img(F)$ are closed under isomorphisms in $\mathcal{C}$.

For a complex $\cpx{X}=(X^i, d_X^i)_{i\in\mathbb{Z}}$ over $\mathcal{C}$, the morphism $d_X^i: X^i\ra X^{i+1}$ is called a
\emph{differential} of $\cpx{X}$. For simplicity, we sometimes write
$(X^i)_{i\in\mathbb{Z}}$ for $\cpx{X}$ without mentioning $d^i_X$. For an integer $n$, we denote by
$\cpx{X}[n]$ the complex by shifting $n$
degrees, namely $(\cpx{X}[n])^i=X^{n+i}$ and $d^i_{\cpx{X}[n]}=(-1)^nd^{n+i}_X$.

Let $\C{\mathcal{C}}$ be the category of all complexes over
$\mathcal{C}$ with chain maps, and $\K{\mathcal{C}}$ the homotopy
category of $\C{\mathcal{C}}$. As usual, we denote by
$\Cf{\mathcal{C}}$ the category of bounded-above complexes over $\mathcal{C}$, and
by $\Kf{\mathcal{C}}$ the homotopy category of $\Cf{\mathcal{C}}$.
Note that $\K{\mathcal{C}}$ and $\Kf{\mathcal{C}}$ are triangulated
categories. For a triangulated category, its shift functor is
denoted by $[1]$ universally.

Now, let $\mathcal{A}$ be an abelian category and $\cpx{X}\in\C{\mathcal{A}}$. For $n\in\mathbb{Z}$, there are two left-truncated complexes
$$
X{^{\leq n}}:\;\; \cdots\lra X^{n-3}\epa{d_X^{n-3}} X^{n-2}\epa{d_X^{n-2}}X^{n-1}\epa{d_X^{n-1}}X^n\lra 0,
$$
$$
\tau^{\leq n}\cpx{X}:\;\; \cdots\lra X^{n-3}\epa{d_X^{n-3}} X^{n-2}\epa{d_X^{n-2}}X^{n-1}\epa{d_X^{n-1}}\Ker(d_X^n)\lra 0,
$$
together with canonical chain maps $\cpx{X}\to X{^{\leq n+1}}\to X{^{\leq n}}$ and $\tau^{\leq n}\cpx{X}\to \tau^{\leq n+1}\cpx{X}\to\cpx{X}$.
Dually, right-truncated complexes $X^{\geq n}$ and $\tau^{\geq n}\cpx{X}$ (by taking cokernels on the left) can be defined.
Further, there are bi-truncated complexes for a pair of integers $(n, m)$ with $n<m$:
$$
X^{[n, m)}:\;\;0\lra X^n\epa{d_X^n} X^{n+1}\lra\cdots \epa{d_X^{m-2}}X^{m-1}\epa{d_X^{m-1}}\Ker(d_X^m)\lra 0,
$$
$$
X^{(n, m]}:\;\;0\lra \Coker(d_X^{n-1})\epa{\overline{d_X^n}} X^{n+1}\lra\cdots \epa{d_X^{m-2}}X^{m-1}\epa{d_X^{m-1}}X^m\lra 0.
$$
where $\overline{d_X^n}$ is induced from $d_X^n$.

Let $H^n(\cpx{X})$ be the
cohomology of $\cpx{X}$ in degree $n$. Then $H^n(-)$ is a functor from $\mathscr{C}(\mathcal{A})$ to $\mathcal{A}$ for all $n$.

Next, we recall the definition of semi-orthogonal decompositions of triangulated categories (for example, see \cite[Chapter 11]{Huybrechts}).
Note that semi-orthogonal decompositions are also termed hereditary torsion pairs (see \cite[Chapter I.2]{BI})
and closely related to Bousfield localizations \cite[Section 9.1]{Neeman}) and $t$-structures of triangulated categories (see \cite{BBD}).

\begin{Def}\label{HTP}
Let $\mathscr{D}$ be a triangulated category. A pair $(\mathscr{X}, \mathscr{Y})$ of full subcategories $\mathscr{X}$ and $\mathscr{Y}$ of $\mathscr{D}$ is called a \emph{semi-orthogonal decomposition} of $\mathscr{D}$ if

$(1)$ $\mathscr{X}$ and $\mathscr{Y}$ are triangulated subcategories of $\mathscr{D}$.

$(2)$ $\Hom_{\mathscr{D}}(X, Y)=0$ for all $X\in\mathscr{X}$ and $Y\in\mathscr{Y}$.

$(3)$ For each object $D\in\mathscr{D}$, there exists a triangle
$$X_D\lraf{f_D} D\lraf{g^D} Y^D\lra X_D[1]$$ in $\mathscr{D}$
such that $X_D\in\mathscr{X}$ and $Y^D\in\mathscr{Y}$.
\end{Def}

The following is well known (for example, see \cite[Chapter I.2]{BI}) and will be used later.

\begin{Lem}\label{Torsion pair}
Let $(\mathscr{X}, \mathscr{Y})$ be a semi-orthogonal decomposition of $\mathscr{D}$. The following statements hold.

$(1)$ The inclusion ${\bf i}: \mathscr{X}\to\mathscr{D}$ has a right adjoint ${\bf R}: \mathscr{D}\to\mathscr{X}$
given by $D\mapsto X_D$ for each $D\in\mathscr{D}$, such that $f_D$ is the counit adjunction morphism of $D$.
Dually, the inclusion ${\bf j}: \mathscr{Y}\to\mathscr{D}$ has a left adjoint ${\bf L}: \mathscr{D}\to\mathscr{Y}$
given by $D\mapsto Y^D$ for each $D\in\mathscr{D}$ such that $g^D$ is the unit adjunction morphism of $D$.

$(2)$ $\Ker({\bf L})=\mathscr{X}$ and ${\bf L}$ induces a triangle equivalence $\overline{\bf L}: \mathscr{D}/\mathscr{X}\to\mathscr{Y}$ of which a
quasi-inverse is the composition of ${\bf j}$ with the localization functor $\mathscr{D}\to\mathscr{D}/\mathscr{X}$.
Dually, $\Ker({\bf R})=\mathscr{Y}$ and ${\bf R}$ induces a triangle equivalence $\overline{\bf R}: \mathscr{D}/\mathscr{Y}\to\mathscr{X}$
of which a quasi-inverse is the composition of ${\bf i}$ with the localization functor $\mathscr{D}\to\mathscr{D}/\mathscr{Y}$.
\end{Lem}

Semi-orthogonal decompositions are intimately connected with recollements of triangulated categories introduced by Beilinson, Bernstein
and Deligne in \cite{BBD} for understanding derived
categories of perverse sheaves over singular spaces. Later, this has widely been used in representation theories of Lie algebras and associative algebras (for example, see \cite{CPS2, happel2, Bz1, XC4}).

\begin{Def}{\rm \cite{BBD}} \label{def01}
Let $\mathcal{D}, \mathcal{D'}$ and $\mathcal{D''}$ be triangulated categories.
The category $\mathcal{D}$ is called a \emph{recollement} of $\mathcal{D'}$ and
$\mathcal{D''}$ (or there is a recollement among $\mathcal{D}'', \mathcal{D}$ and
$\mathcal{D}'$) if there are six triangle functors
$$\xymatrix{\mathcal{D''}\ar^-{i_*=i_!}[r]&\mathcal{D}\ar^-{j^!=j^*}[r]
\ar^-{i^!}@/^1.2pc/[l]\ar_-{i^*}@/_1.6pc/[l]
&\mathcal{D'}\ar^-{j_*}@/^1.2pc/[l]\ar_-{j_!}@/_1.6pc/[l]}$$ satisfying the four properties:

$(1)$ The four pairs $(i^*,i_*),(i_!,i^!),(j_!,j^!)$ and $(j^*,j_*)$ of functors are adjoint.

$(2)$ The three functors $i_*,j_*$ and $j_!$ are fully faithful.

$(3)$ $i^!j_*=0$ (and thus also $j^! i_!=0$ and $i^*j_!=0$).

$(4)$ There are two triangles for each object $D$  in
$\mathcal D$:
$$
i_!i^!(D)\lra D\lra j_*j^*(D)\lra i_!i^!(D)[1],
$$
$$
j_!j^!(D)\lra D\lra i_*i^*(D)\lra j_!j^!(D)[1],
$$
where $i_!i^!(D)\ra D$ and $j_!j^!(D)\ra D$ are counit adjunction morphisms, and where $D\ra j_*j^*(D)$ and
$D\ra i_*i^*(D)$ are unit adjunction morphisms.
\end{Def}

If $\mathcal{D}$ is a recollement of $\mathcal{D'}$ and
$\mathcal{D''}$, then the pairs $\big(\Img(j_!),\Img(i_*)\big)$ and $\big(\Img(i_*),\Img(j_*)\big)$ are semi-orthogonal decomositions of $\mathscr{D}$.

Next, we recall the definition of homotopy colimits and limits in triangulated categories.

\begin{Def}{\rm \cite{BN, Neeman}}\label{def-holim}
Let $\mathscr{D}$ be a triangulated category such that countable coproducts exist in $\mathscr{D}$.  Let
$$X_0\lraf{f_0}X_1\lraf{f_1}X_2\lraf{f_2}\cdots\lraf{f_n}X_{n+1}\lra\cdots$$
be a sequence of objects and morphisms in $\mathscr{D}$. The \emph{homotopy colimit} of this sequence, denoted by $\underrightarrow{\hocolim}(X_n)$,
is given, up to non-canonical isomorphism, by the triangle
$$
\bigoplus_{n\geq 1}X_n\lraf{(1-f_*)} \bigoplus_{n\geq 1}X_n\lra \underrightarrow{\hocolim}(X_n)\lra \bigoplus_{n\geq 1}X_n[1]
$$
where the morphism $(1-f_*)$  is induced by  $({\rm Id}_{X_i}, -f_i): X_i\to X_i\oplus X_{i+1}\subseteq \bigoplus_{n\geq 1}X_n$ for all $i\in\mathbb{N}$.

Dually, the \emph{homotopy limit}, denoted by $\underleftarrow{\holim}$, of a sequence of objects and morphisms in a triangulated category with countable products can be defined.
\end{Def}

Now, we consider homotopy colimits and limits in derived categories of abelian categories.
Let $\mathcal{A}$ be an abelian category. Recall that $\mathcal{A}$ \emph{satisfies $\rm{AB4}$} if coproducts indexed over sets exist in $\mathcal{A}$ and
coproducts of short exact sequences in $\mathcal{A}$ are exact. Dually, $\mathcal{A}$ \emph{satisfies} $\rm{AB4'}$ if products indexed over sets exist in $\mathcal{A}$ and products of short exact sequences in $\mathcal{A}$ are exact. A typical example of abelian categories satisfying both $\rm{AB4}$ and $\rm{AB4'}$ is the module category of a ring.

A full subcategory $\mathcal{B}$ of the abelian category $\mathcal{A}$ is called an \emph{abelian subcategory} of $\mathcal{A}$ if $\mathcal{B}$ is an abelian category and the inclusion $\mathcal{B}\to\mathcal{A}$ is an exact functor between abelian categories. This is equivalent to saying that $\mathcal{B}$ is closed under taking
kernels and cokernels in $\mathcal{A}$.

The following result is taken from \cite[Lemma 6.1]{Neeman1} and its dual statement.

\begin{Lem}\label{Homotopy}
Let $\mathcal{A}$ be an abelian category satisfying $AB4$ and $AB4'$.
For a complex $\cpx{X}\in\C{\mathcal{A}}$, there are isomorphisms in $\K{\mathcal{A}}$ (and also in $\D{\mathcal{A}}$)
$$
\cpx{X}\simeq \underrightarrow{\hocolim}(X^{\geq -n})\simeq \underrightarrow{\hocolim}(X^{[-n, n+1)})\simeq\underleftarrow{\holim}(X^{\leq n})\simeq \underleftarrow{\holim}(X^{(-n, n+1]}),
$$where $n$ runs over all natural numbers.
\end{Lem}

Finally, we mention a general fact on adjoint pairs of functors.

\begin{Lem}\label{fact-adj} Let $F: \mathcal{C}\ra\mathcal{D}$ and $G:\mathcal{D}\ra\mathcal{C}$ be two functors of categories such that $(F,G)$ forms a adjoint pair.

$(1)$ If $G$  is fully faithful, then the counit $\epsilon: F\circ G\ra Id_{\mathcal{D}}$ is an isomorphism.

$(2)$ $G$ gives rise to an equivalence $\mathcal{D}\ra \Img(G)$ with the quasi-inverse given by the restriction of $F$ to $\Img(G)$. In particular, if $C\in \Img(G)$, then the unit $\eta_C: C\ra GF(C)$ is an isomorphism.
\end{Lem}

\subsection{Derived categories of exact categories}\label{Exact category}

An \emph{exact category} (in the sense of Quillen) is by definition an additive category endowed with a class of conflations
closed under isomorphism and satisfying certain axioms (for example, see \cite[Section 4]{Keller}). When the additive category is an abelian
category, the class of conflations coincides with the class of short exact sequences.

Let $\mathcal{E}$ be an exact category and $\mathcal{F}$ a full subcategory of $\mathcal{E}$.
Suppose that $\mathcal{F}$ is \emph{closed under extensions} in $\mathcal{E}$,
that is, for any conflation $0\to X\to Y\to Z\to 0$ in $\mathcal{E}$ with both $X,Z\in \mathcal{F}$, we have $Y\in\mathcal{F}$. Then $\mathcal{F}$, endowed with the conflations in $\mathcal{E}$ having their terms in $\mathcal{F}$, is an exact
category, and the inclusion $\mathcal{F}\subseteq\mathcal{E}$ is a fully faithful exact functor. In this case,
$\mathcal{F}$ is called a \emph{fully exact subcategory} of $\mathcal{E}$ (see \cite[Section 4]{Keller}).

From now on, let $\mathcal{A}$ be an abelian category and $\mathcal{E}$ a fully exact subcategory of $\mathcal{A}$.  A complex $\cpx{X}\in\C{\mathcal{E}}$ is said to be \emph{strictly exact}
if it is exact in $\C{\mathcal{A}}$ and all of its boundaries belong to $\mathcal{E}$.
Let $\mathscr{K}_{\rm {ac}}(\mathcal{E})$ be the full subcategory of $\K{\mathcal{E}}$ consisting of those complexes which are isomorphic to
strictly exact complexes. Then $\mathscr{K}_{\rm {ac}}(\mathcal{E})$ is a full triangulated subcategory of $\K{\mathcal{E}}$ closed under direct summands.
The \emph{unbounded derived category} of $\mathcal{E}$, denoted by $\D{\mathcal{E}}$, is defined to be the Verdier quotient of $\K{\mathcal{E}}$
by $\mathscr{K}_{\rm {ac}}(\mathcal{E})$. Similarly, the bounded-below, bounded-above and bounded derived categories
$\Dz{\mathcal{E}}$, $\Df{\mathcal{E}}$ and $\Db{\mathcal{E}}$ can be defined. Observe that the canonical functor $\mathscr{D}^*(\mathcal{E})\to\D{\mathcal{E}}$ is fully faithful
for $*\in\{+, -, b\}$.

If $F: \mathcal{E}_1\to \mathcal{E}_2$ is an additive functor of exact categories, then $F$ induces a functor $K(F): \mathscr{K}(\mathcal{E}_1)\to \mathscr{K}(\mathcal{E}_2)$ of homotopy categories. Further, if $F$ is an exact functor, that is, $F$ sends conflations in $\mathcal{E}_1$ to the ones in $\mathcal{E}_2$, then $F$ induces a functor $D(F): \mathscr{D}(\mathcal{E}_1)\to \mathscr{D}(\mathcal{E}_2)$ of derived categories.

Let $\mathcal{F}$ and $\mathcal{E}$ be fully exact subcategories of $\mathcal{A}$ with
$\mathcal{F}\subseteq \mathcal{E}$. Then $\mathcal{F}$ can be regarded as a fully exact subcategory of $\mathcal{E}$.
We consider the following three conditions:

$(a)$ If $0\to X\to Y\to Z\to 0$ is an exact sequence in $\mathcal{A}$ with $X\in\mathcal{E}$ and $Y, Z\in \mathcal{F}$, then $X\in\mathcal{F}$.

$(b)$ Each exact sequence $0\to E_1\to E_0\to F\to 0$ in $\mathcal{A}$ with $E_1, E_0\in\mathcal{E}$ and $F\in\mathcal{F}$ can be completed into
an exact commutative diagram
$$
\xymatrix{0\ar[r]&F_1\ar[r]\ar[d]&F_0\ar[r]\ar[d]& F\ar[r]\ar@{=}[d]&0\\
0\ar[r]&E_1\ar[r]& E_0\ar[r]& F\ar[r]&0
}
$$
with $F_1, F_0\in\mathcal{F}$.

$(c)$ There is a natural number $n$ such that, for each object $E\in\mathcal{E}$, there is a long exact sequence in $\mathcal{A}$:
$$0\lra F_n\lraf{f_n} \cdots\lraf{f_2} F_1\lraf{f_1} F_0\lraf{f_0} E\lra 0$$
with $F_i\in\mathcal{F}$ and $\Img(f_i)\in\mathcal{E}$ for all $0\leq i\leq n$.

The following result describes impacts of these conditions on derived categories.
\begin{Lem}\label{resolution}
$(1)$ If $(a)$ and $(b)$ hold, then the inclusion $\mathcal{F}\subseteq\mathcal{E}$ induces a fully faithful triangle functor $\Df{\mathcal{F}}\to\Df{\mathcal{E}}$.

$(2)$ Suppose that $\mathcal{E}$ is closed under direct summands in $\mathcal{A}$. If $(a)$ and $(c)$ hold, then the inclusion $\mathcal{F}\subseteq\mathcal{E}$ induces a triangle equivalence $\D{\mathcal{F}}\ra \D{\mathcal{E}}$ which restricts to an equivalence $\mathscr{D}^*(\mathcal{F})\ra \mathscr{D}^*(\mathcal{E})$ for any $*\in\{+, -, b\}$.
\end{Lem}

{\it Proof.} $(1)$ follows from the dual of \cite[Theorem 12.1]{Keller} (see also \cite[Proposition A.2.1]{Pos1}), while $(2)$ follows from
\cite[Proposition A.5.6]{Pos1}. $\square$

If all arrows in $(a), (b)$ and $(c)$ are reversed, we get  dual conditions of $(a), (b)$ and $(c)$, respectively. So the dual version of Lemma \ref{resolution} holds true.

\subsection{Derived functors of the module categories of rings}

Let $R$ be a (unitary associative) ring. The full subcategories of projective and injective $R$-modules
are denoted by $\Pmodcat{R}$ and $\Imodcat{R}$, respectively.
If $M$ is an $R$-module and $I$ is a nonempty set, then $M^{(I)}$ and $M^I$ denote the direct sum and product of $I$ copies of $M$,
respectively. The projective dimension and the endomorphism ring of $M$ are denoted by $\pd({_R}M)$ and $\End_R(M)$, respectively.

A full subcategory $\mathcal{T}$ of $R\Modcat$ is called a \emph {thick subcategory} if it is closed under direct summands in $R\Modcat$ and has
the \emph {two out of three property}: for any short exact sequence $0\to X\to Y\to Z\to 0$ in $R\Modcat$ with two terms in $\mathcal{T}$, the third term belongs to $\mathcal{T}$ as well.

In the sequel, we simply write $\C{R}$, $\K{R}$ and $\D{R}$ for
$\C{R\Modcat}$, $\K{R\Modcat}$ and $\D{R\Modcat}$, respectively, and
regard $R\Modcat$ as the subcategory of $\D{R}$ consisting of
all stalk complexes concentrated in degree zero.

Let $\K{R}_P$ (respectively, $\K{R}_I$) be the smallest
full triangulated subcategory of $\K{R}$, which
contains all the bounded-above (respectively, bounded-below)
complexes of projective (respectively, injective) $R$-modules, and
is closed under arbitrary direct sums (respectively, direct
products).
Since $\K{\Pmodcat{R}}$ satisfies theses two properties, we have $\K{R}_P\subseteq\K{\Pmodcat{R}}$.
Let $\mathscr{K}_{\rm {ac}}(R)$ be the full subcategory of $\K{R}$ consisting of all
exact complexes. Then $(\K{R}{_P}, \mathscr{K}_{\rm {ac}}(R))$ forms a semi-orthogonal decomposition
of $\K{R}$. Let $_p(-):\K{R}\to\K{R}{_P}$ be a right adjoint of the inclusion
$\K{R}_P\to\K{R}$. Then, by Lemma \ref{Torsion pair}, the functor $_p(-)$ induces a triangle equivalence $\D{R}\lraf{\simeq}\K{R}{_P}$ of which
a quasi-inverse is the composition of the inclusion $\K{R}_P\to\K{R}$ with the localization functor $\K{R}\to \D{R}$.
Moreover, for each $\cpx{X}$ in $\K R$, the associated counit adjunction morphism $\alpha_{\cpx X}:{_p}\cpx{X}\to \cpx{X}$ is a quasi-isomorphism. Thus
$\alpha_{\cpx X}$ or simply $_p\cpx{X}$ is called a \emph{projective resolution} of $\cpx{X}$ in $\D{R}$. For example, if $X$ is  an
$R$-module, then $_pX$ can be chosen as a deleted projective resolution of $_RX$.
Dually, $(\mathscr{K}_{\rm {ac}}(R), \K{R}{_I})$ is a semi-orthogonal decomposition of $\K{R}$. In particular, there exists a quasi-isomorphism $\beta_{\cpx{X}}: \cpx{X}\to{_i}\cpx{X}$ in $\K{R}$ with ${_i}\cpx{X}\in \K{R}_I$. The
complex $_i\cpx{X}$ is called the \emph{injective resolution} of
$\cpx{X}$ in $\D{R}$. In particular, $\Hom_{\K{R}}(\cpx{P},\cpx{X})\simeq \Hom_{\D{R}}(\cpx{P},\cpx{X})$ and $\Hom_{\K{R}}(\cpx{X},\cpx{I})\simeq \Hom_{\D{R}}(\cpx{X}, \cpx{I})$ for $\cpx{P}\in \mathscr{K}^-(\Pmodcat{R}), \cpx{I}\in \mathscr{K}^+(R\mbox{-Inj})$ and $\cpx{X}\in \K{R}.$

Now, let $S$ be another ring. For a triangle functor $F:\K{R}\to\K{S}$, its \emph{total left-derived functor} ${\mathbb L}F:\D{R}\to\D{S}$ is defined by
$\cpx{X}\mapsto F(_p\cpx{X})$, and its \emph{total right-derived
functor} ${\mathbb R}F:\D{R}\to\D{S}$ is defined by $\cpx{X}\mapsto
F(_i\cpx{X})$. Further, if $F(\cpx{X})$ is exact whenever $\cpx{X}$ is exact, then $F$
induces a triangle functor $D(F):\D{R}\to\D{S}$,
$\cpx{X}\mapsto F(\cpx{X})$. In this case, up to natural
isomorphism, ${\mathbb L}F={\mathbb R}F=D(F)$, and
$D(F)$ is called the \emph{derived functor} of $F$.

Given a complex $\cpx{M}$ of $R$-$S$-bimodules, we denote by ${\mathbb R}\Hom_R(\cpx{M},-)$ and $\cpx{M}\otimesL_S-$ the
total right-derived functor of $\cpx{\Hom}_R(\cpx{M},-)$ and the total left-derived functor of
$\cpx{M}\cpx{\otimes}_S-$, respectively. Then $\big(\cpx{M}\otimesL_S-, {\mathbb R}\Hom_R(\cpx{M},-)\big)$ is an adjoint pair of triangle functors. In case of
$Y\in S\Modcat$  and $X\in R\Modcat$, we write $\cpx{M}\otimes_SY$ and $\Hom_R(\cpx{M},\,X)$
for $\cpx{M}\cpx{\otimes}_SY$ and $\cpx{\Hom}_R(\cpx{M},\,X)$, respectively.

\section{Symmetric subcategories of abelian categories\label{sect3}}

Now we introduce $n$-symmetric subcategories of an ableian category for $n\ge 0$.

Recall that an abelian category is \emph{complete} (respectively, \emph{cocomplete}) if it has products (respectively, coproducts) indexed over sets; and \emph{bicomplete} if it is complete and cocomplete.

\begin{Def} \label{sym-cat} Let $n\in \mathbb{N}$, and let $\mathcal{A}$ be a bicomplete abelian category.  An additive full subcategory $\mathcal{B}$ of $\mathcal{A}$ is said to be \emph{$n$-symmetric} if

$(1)$ $\mathcal{B}$ is closed under extensions, products and coproducts.

$(2)$ Given an exact sequence $0\ra X\ra M_{n}\ra\cdots\ra M_1\ra M_0\ra Y\ra 0$ in $\mathcal{A}$, we have $X,Y\in \mathcal{B}$ whenever all $M_i\in \mathcal{B}$.
\end{Def}

\begin{Rem}{\rm Let $\mathcal{B}$ be an additive full subcategory of a bicomplete abelian category $\mathcal{A}$.

$(1)$ If $\mathcal{B}$ is an $n$-symmetric subcategory of $\mathcal{A}$, then $\mathcal{B}$ is an exact, thick subcategory of $\mathcal{A}$. It is also $(n+1)$-symmetric.

$(2)$ If $\mathcal{B}_i$ is an $m_i$-symmetric subcategories of  $\mathcal{A}$ for $i=1,2$ then $\mathcal{B}_1\cap\mathcal{B}_2$ is a max$\{m_1,m_2\}$-symmetric subcategory of $\mathcal{A}$.

$(3)$ If $\mathcal{B}$ is closed under extensions and satisfies Definition \ref{sym-cat}(2), then $\mathcal{B}$ is an $n$-wide subcategory of $\mathcal{A}$ in the sense of \cite[Definition 4.1]{mntty}.

$(4)$ $\mathcal{B}$ is $0$-symmetric if and only if $\mathcal{B}$ is a Serre subcategory (that is, closed under subobjects, quotient objects and extensions) closed under coproducts and products if and only $\mathcal{B}$ is a localizing subcategory (that is, a Serre subcategory closed under coproducts) closed under products.
$\mathcal{B}$ is $1$-symmetric if and only if
$\mathcal{B}$ is an abelian subcategory closed under extensions, coproducts and products.

In particular, if $R$ is  a ring and $\mathcal{A}=R\Modcat$, then a full subcategory $\mathcal{B}$ of $\mathcal{A}$ is $0$-symmetric if and only if there  is an ideal $I$ of $R$ with $I^2 = I$ such that $\mathcal{B}= (R/I)\Modcat$ = $\{M\in R\Modcat \mid IM=0\}$. This follows from \cite[Proposition 6.12]{stenstroem}. Moreover, a full subcategory $\mathcal{B}$ of $\mathcal{A}$ is $1$-symmetric if and only if there is a ring epimorphism $\lambda: R\ra S$ with $\Tor_1^R(S,S)=0$ such that $\mathcal{B}$ coincides with the image of the induced fully faithful functor $\lambda_*: S\Modcat \ra R\Modcat $, that is, $\mathcal{B}=\Img(\lambda_*)\simeq S\Modcat.$ This can be seen by \cite[Lemma 2.1]{XC1}.
}\end{Rem}

\begin{Prop} Let $\mathcal{A}$ and $\mathcal{C}$ be bicomplete abelian categories satisfying AB4 and AB4$'$, and let $F: \D{\mathcal{A}}\ra \D{\mathcal{C}}$ be a triangle functor commuting with coproducts and products. If there are integers $s\le r$ such that
$H^i(FX)=0$ for all $X\in \mathcal{A}$, $r< i$, or $i < s$, then $\mathscr{E}:= \mathcal{A}\cap\Ker(F)$ is an $(r-s)$-symmetric subcategory of $\mathcal{A}$.
\label{sym-pro1}
\end{Prop}

{\it Proof.} It is easy to see that $\mathscr{E}$ is an additive subcategory and closed under extensions, coproducts and products. Let $n:=r-s$ and let $0\ra X\ra M_n\ra \cdots\ra M_1\ra M_0\ra Y\ra 0$ be an exact sequence in $\mathcal{A}$ with $M_j\in \mathscr{E}$ for all $j$. By assumption, we have $ H^i(F(Y))=0$ for $i>r$ or $i<s$. Now, let $s\le i\le r$. Then $n+1+i\ge n+1+s=r+1$. It follows from $F(M_j)=0$ for $0\le j\le n$ that $FY\simeq F(X)[n+1]$ in $\D{\mathcal{C}}.$ This implies
$ H^i(F(Y))\simeq H^i(F(X)[n+1])\simeq H^{n+1+i}(F(X))$ for $i\in\mathbb{Z}$.
In particular, $H^i(F(Y))=0$ by assumption. Thus $H^i(F(Y))=0$ for all $i\in \mathbb{Z}$ and $Y\in \Ker(F)$. This also shows $X\in \Ker(F)$. Thus $\mathscr{E}$ is an $n$-symmetric subcategory of $\mathcal{A}$. $\square$

As an application of Proposition \ref{sym-pro1}, we consider the module categories of rings.

\begin{Bsp} \label{bsp3.4}
{\rm Let $R$ and $S$ be rings, and let ${_R}M_S$ be an $R$-$S$-bimodule.

(1) If $M_S$ has a finite projective resolution of length $n$ by finitely generated projective right $S$-modules, that is, there is an extact sequence $0\ra Q_n\ra \cdots\ra Q_1\ra Q_0\ra M_S\ra 0$ with all $Q_j$ being finitely generated projective right $S$-modules, then
$\Ker(M\otimesL_S-)\cap S\Modcat$ is an $n$-symmetric subcategory of $S\Modcat$.

(2) If ${_R}M$  has a finite projective resolution of length $n$ by finitely generated projective modules, then $\Ker(\mathbb{R}\Hom_R(M,-))\cap R\Modcat$ is an $n$-symmetric subcategory of $R\Modcat$.

Note that $\Ker(M\otimesL_S-)\cap S\Modcat$ = $\{Y\in S\Modcat\mid \Tor_i^S(M,Y)=0, \forall i\ge 0 \}$ and
$\Ker(\mathbb{R}\Hom_R(M,-))\cap R\Modcat$ = $\{X\in R\Modcat\mid \Ext^i_R(M,X)=0, \forall i\ge 0 \}$.

(3) Let $_AT$ be a good $n$-tilting $A$-module with $B :=\End_A(T)$. Then $\mathscr{E}:=\Ker(T\otimesL_B-)\cap{B\Modcat}$ is always an $n$-symmetric subcategory of $B\Modcat$. It follows from \cite[Theorem 5.12 and Theorem 7.5]{Bz1} that $\mathscr{E}$ is $0$-symmetric if and only if $_AT$ is pure-projective if and only if the heart of the $t$-structure induced from the tilting module $_AT$ is a Grothendieck category. Furthermore, it follows from \cite[Theorem 1.1]{XC2} that $\mathscr{E}$ is $1$-symmetric if and only if $_AT$ is a homological tilting module.
}\end{Bsp}
A further discussion on $n$-symmetric subcategories will be done in a forthcoming article.

\section{Tilting modules and derived recollements\label{sect4}}

This entire section is devoted to a proof of Theorem \ref{Main result}. We first make a couple of preparations

\subsection{Derived functors induced by good tilting modules\label{sect4.1}}

Throughout this section, $A$ denotes a unitary ring, $T$ a {\bf{good}} $n$-tilting $A$-module with $(T1), (T2)$, and $(T3)$ for a natural
number $n$, and $B:=\End_A(T)$. Further, we define
$$
G:={_A}T\otimesL_B- : \,\D{B}\lra \D{A}\quad\mbox{and}\quad H:=\rHom_A(T,-):\,
\D{A}\lra \D{B},
$$
to be the total left- and right-derived functors of ${_A}T{_B}$, respectively.
We write by $$\eta': Id_{B\Modcat}\lra \Hom_A(T,T\otimes_B-)\quad\mbox{and}\quad \eta:Id_{\D{B}}\lra H\circ G $$
the unit adjunctions associated with the adjoint pairs $(T\otimes_B-, \Hom_A(T,-))$ and $(G, H)$, respectively.
Recall that the restriction of the localization functor $\K{B}\to\D{B}$ to $\K{B}_p$ is a triangle equivalence.
A quasi-inverse of this equivalence is denoted by $$_p(-): \D{B}\lra \K{B}_P.$$ Without loss of generality,
we assume that ${_p}M$ is a deleted projective resolution of $M$ for each $B$-module $M$.

The derived categories $\D{A}$ and $\D{B}$ are related by the following recollement (see \cite[Theorem 2.2]{Bz2}). This is used in \cite[Lemma 5.2]{XC2} to understand when a tilting module is homological.

\begin{Theo}{\rm \cite{Bz2}}\label{Torsion}
There exists a recollement of triangulated categories:
$$
\xymatrix@C=1.2cm{\Ker(G) \ar[r]
&\D{B}\ar[r]^{G}\ar@/^1.2pc/[l]\ar@/_1.2pc/[l]
&\D{A}\ar@/^1.2pc/[l]^H\ar@/_1.2pc/[l]}.$$
In particular, $(\Ker(G), \Img(H))$ is a semi-orthogonal decomposition of $\D{B}$.
\end{Theo}

\medskip
Observe that the modules ${_A}T$ and $T_B$ have finite projective dimension. In fact, $T_B$ has a finite projective resolution by finitely generated projective right $B$-modules since the tilting module $_A T$ is good. The next result is deduced from \cite[Theorem 10.5.9, Corollary 10.5.11]{We}.

\begin{Lem}\label{HT}
$(1)$ Let $\alpha_{\cpx{Y}}: {_p}\cpx{Y}\to \cpx{Y}$ be the projective resolution of $\cpx{Y}$ in $\D{B}$.  If $\Tor^B_j(T, Y^i)=0$ for all $i\in\mathbb{Z}$ and $j\geq 1$, then $T\otimes_B\alpha_{\cpx{Y}}: G(\cpx{Y})=T\otimes_B{_p}\cpx{Y}\ra T\otimes_B\cpx{Y}$ is an isomorphism in $\D{A}$.

\smallskip
$(2)$ Let $\beta_{\cpx{X}}: \cpx{X}\to {_i}\cpx{X}$ be the injective resolution of $\cpx{X}$ in $\D{A}$. If $\Ext_A^j(T, X^i)=0$ for all $i\in\mathbb{Z}$ and $j\geq 1$, then $\Hom_A(T, \beta_{\cpx{X}}): \Hom_A(T, \cpx{X})\ra \Hom_A(T, {_i}\cpx{X})=H(\cpx{X})$ is an isomorphism in $\D{B}$.
\end{Lem}

Combining Theorem \ref{Torsion} with Lemma \ref{HT}, we have the lemma.

\begin{Lem}\label{1-tilting}
$(1)$ For each $P\in\Pmodcat{B}$, the unit adjunction morphism $\eta_P': P \to \Hom_A(T, T\otimes_BP)$ is
injective with $\Coker(\eta_P')\in \mathscr{E}$.

$(2)$ For each $\cpx{M}\in\D{B}$, there is a commutative
diagram of triangles in $\D{B}$:
$$
\xymatrix{
\Coker(\eta'_{{_p}\cpx{M}})[-1]\ar[r]\ar@{=}[d]&{_p}\cpx{M}\ar[r]^-{\eta'_{{_p}\cpx{M}}}\ar[d]_-{\simeq}^-{\alpha_{\cpx{M}}}
&\Hom_A(T, T\otimes_B {_p}\cpx{M})\ar[d]_-{\simeq}^-{\Hom_A(T, \beta_{T\otimes_B{_p}\cpx{M}})}\ar[r]& \Coker(\eta'_{{_p}\cpx{M}})\ar@{=}[d]\\
\Coker(\eta'_{{_p}\cpx{M}})[-1]\ar[r]&\cpx{M}\ar[r]^-{\eta_M} & H\circ G(\cpx{M})\ar[r]&\Coker(\eta'_{{_p}\cpx{M}})}$$
where $\eta'_{{_p}\cpx{M}}:=(\eta_{{_p}M^n}')_{n\in\mathbb{Z}}$ and $\Coker(\eta'_{{_p}\cpx{M}})\in\Ker(G)$.

\smallskip
$(3)$ $$\Ker(G)=\{\cpx{\overline{Y}}\in\D{B}\mid\cpx{\overline{Y}}\simeq \cpx{Y} {\rm \; in \;}\D{B}
{\rm \; with \; } Y^i\in\mathscr{E} {\rm \; for \ all \;}
i\in\mathbb{Z}\}\quad \mbox{and}
$$
$$\Img(H)=\{\cpx{\overline{Z}}\in\D{B}\mid\cpx{\overline{Z}}\simeq
\cpx{Z} \mbox{\; in \;}\D{B} \mbox{\; with\;}
Z^i\in\Hom_A(T,\Add(T))\mbox{\;for\, all \,} i\in\mathbb{Z}\},$$
where $\Hom_A(T,\Add(T)):=\{\Hom_A(T, T')\in B\Modcat\mid T'\in \Add({_A}T)\}$.
\end{Lem}

{\it Proof.} $(1)$ Since the canonical map
$\Hom_A(T, T)^{(\alpha)}\to \Hom_A(T, T^{(\alpha)})$ is injective for any nonempty set $\alpha$, the map $\eta_P'$ is injective for any free $B$-module $P$, and therefore $\eta_P'$ is injective for any projective $B$-module $P$. By Lemma \ref{HT},
$H\circ G(P)\simeq \Hom_A(T, T\otimes_BP)$. It follows from Theorem \ref{Torsion} that there exists a triangle
$$\cpx{X}_P \lra P\lraf{\eta_P'} \Hom_A(T, T\otimes_BP)\lra\cpx{X}_P [1]$$ in $\D{B}$ with $\cpx{X}_P\in\Ker(G)$.
As $\eta_P'$ is injective, we have $\Coker(\eta_P')\simeq \cpx{X}_P[1]\in\Ker(G)$, and therefore $\Coker(\eta_P')\in\mathscr{E}$.

$(2)$ This follows from Theorem \ref{Torsion}, Definition \ref{def01}(4) and Lemma \ref{HT}(2) because ${_p}\cpx{M}\in\K{B}_P\subseteq \K{\Pmodcat{B}}$ and
$\Ext_A^j(T, T^{(\alpha)})=0$ for $j\geq 1$ and nonempty sets $\alpha$.

$(3)$ This is shown in {\rm \cite[Proposition 4.6]{XC1}} for good $1$-tilting modules, but the proof
there also works for good $n$-tilting modules. $\square$

\medskip
Thanks to \cite[Theorem 3.5]{Bz1}, there exists a t-structure in $\D{A}$ associated with ${_A}T$ such that its heart is given by
$$\mathcal{H}:=\{X\in \D{A}\mid \Hom_{\D{A}}(T, X[n])=0\mbox{\;for\, all \,} n\neq 0\}.$$
Then $\mathcal{H}$ contains $\Add(_AT)$ and is an abelian category of which projective objects are
isomorphic in $\D{A}$ to modules in $\Add({_A}T)$.

\begin{Lem}{\rm \cite[Proposition 5.5]{Bz1}} \label{Baz}
The following are true.

$(1)$ The restriction of the functor $H$ to $\mathcal{H}$ yields an exact and fully faithful functor
$$H_T=\Hom_\mathcal{H}(T,-): \mathcal{H}\lra B\Modcat.$$

$(2)$ $\Img(H_T)=\Img(H)\cap B\Modcat$.

$(3)$ $H_T$ has a left adjoint $F: B\Modcat\to \mathcal{H}$ given by the composition of the restriction of $G$
to $B\Modcat$ with the left adjoint of the inclusion $\mathcal{H}\to \mathscr{U}$, where
$\mathscr{U}:=\{X\in \D{A}\mid \Hom_{\D{A}}(T, X[n])=0\mbox{\;for\,}n>0\}.$
\end{Lem}

Let $\mathcal{Y}:=\Img(H_T)$. Then $H_T$ induces an equivalence $\mathcal{H}\simeq \mathcal{Y}$. We have the following.

\begin{Koro}\label{Restriction}
$(1)$ The restriction of $H_T$ to $\Add({_A}T)$ coincides with the one of $\Hom_A(T,-)$ to $\Add(_AT)$,  and the restriction of $F$, defined in Lemma \ref{Baz}(3), to $\Pmodcat{B}$ coincides with the one of $T\otimes_B-$ to $\Pmodcat{B}$.

$(2)$ $\mathcal{Y}$ is an abelian subcategory of $B\Modcat$ closed under isomorphisms, extensions and direct products.
Moreover, a $B$-module $M$ is a projective object of $\mathcal{Y}$ if and only if $M\simeq\Hom_A(T, T')$ for some $T'\in\Add({_A}T)$.
\end{Koro}

{\it Proof.} (1) If $P\in\Pmodcat{B}$, then $G(P)=T\otimes_BP\in\Add({_A}T)\subseteq\mathscr{U}$. Now, Corollary \ref{Restriction} follows from Lemma \ref{Baz}(1) and $(3)$.

(2) By Lemma \ref{Baz}(1), we have the equivalence $H_T: \mathcal{H}\ra \mathcal{Y}$, while the projective objects of $\mathcal{H}$ are just objects of $\Add(_AT)$, up to ismorphism in $\D{A}$. Hence (2) follows. $\square$

\medskip
Now, let
$$\mathscr{E}:=\{Y\in B\Modcat \mid \Tor_i^B(T,Y)=0\mbox{\;for\, all \,} i\geq0\}=B\Modcat\cap\Ker(G).$$
In general, $\mathscr{E}$ is not an abelian subcategory of $B\Modcat$. Nevertheless, we have the property.

\begin{Lem}\label{thick}
$(1)$ $\mathscr{E}$ is an $n$-symmetric subcategory of $B\Modcat$. Particularly, $\mathscr{E}$ is a thick subcategory of $B\Modcat$.

$(2)$ If a $B$-module $M$ is quasi-isomorphic in $\C{B}$ to a complex with all terms in $\mathscr{E}$, then $M$ belongs to $\mathscr{E}$.

$(3)$ Let $\mathscr{E}_0:=\{\Coker(\eta'_P)\mid P\in\Pmodcat{B}\}$. Then a $B$-module $M\in\mathscr{E}$ if and only if there exists an exact sequence
$0\to M\to E_\infty\to E^1\to 0$ of $B$-modules such that  $E^1\in\mathscr{E}_0$ and $E_\infty$ admits a long exact sequence
$$
\cdots\lra E_i\lra\cdots\lra E_2\lra E_1\lra E_\infty\lra 0
$$
with $E_i\in\mathscr{E}_0$ for all $i\ge 1$.
\end{Lem}

{\it Proof.} $(1)$ follows from Example \ref{bsp3.4}(3).
$(2)$ is due to the characterization of $\Ker(G)$ in
Lemma \ref{1-tilting}(3).

(3) The sufficiency of $(3)$ follows from $(1)$. To show the necessity of $(3)$,  we take $M\in\mathscr{E}$. Then $G(M)=0$.
Let $_pM:=(P^{-i})_{i\in\mathbb{N}}$ be a deleted projective resolution of $M$.
By Lemma \ref{1-tilting}(2), $M[1]\simeq \Coker(\eta'_{{_p}M})$ in $\D{B}$.
Define $E_i:=\Coker(\eta'_{P^{-i}})$ for each $i\in\mathbb{N}$. Then $\Coker(\eta'_{{_p}M})$ has the form $\cdots \to E_i\to E_{i-1}\to\cdots\to E_2\lraf{d_2} E_1\lraf{d_1} E_0\to 0$, where $E_i$ is of the degree $(-i)$. Since $P^{-i}\in\Pmodcat{B}$, we have $E_i\in\mathscr{E}_0$ by Lemma \ref{1-tilting}(1).
Note that $H^i(\Coker(\eta'_{{_p}M}))\simeq H^{i+1}(M)=0$ for any $i\neq -1$ and $H^{-1}(\Coker(\eta'_{{_p}M}))\simeq H^0(M)=M$.
This implies that $d_1$ is surjective and the kernel of the map $\Coker(d_2)\to E_0$ induced from $d_1$ is isomorphic to $M$.
Now, if $E_\infty:=\Coker(d_2)$ and $E^1:=E_0$, then the necessity of $(3)$ holds. $\square$

\medskip
The category $\mathscr{E}$ can be regarded as a fully exact subcategory of $B\Modcat$. Moreover,
a complex $\cpx{X}\in\C{\mathscr{E}}$ is strictly exact (see Section \ref{Exact category}) if and only if it is exact in $\C{B}$ by Lemma \ref{thick}(2).
This implies that $\mathscr{K}_{\rm {ac}}(\mathscr{E})$ is the full triangulated subcategory of $\K{\mathscr{E}}$ consisting of exact complexes.
Moreover, by Lemma \ref{thick}(1), $\mathscr{K}_{\rm {ac}}(\mathscr{E})$ is closed under
arbitrary direct sums and products in $\K{B}$, and therefore in $\K{\mathscr{E}}$. Recall that
the \emph{unbounded derived category} $\D{\mathscr{E}}$ of $\mathscr{E}$ is defined to be the Verdier quotient of $\K{\mathscr{E}}$
by $\mathscr{K}_{\rm {ac}}(\mathscr{E})$. Then $\D{\mathscr{E}}$ has direct sums and products, and the localization functor $\K{\mathscr{E}}\to\D{\mathscr{E}}$
preserves direct sums and products by \cite[Lemma 1.5]{BN} and its dual statement.

Let $$i:\mathscr{E}\lra B\Modcat\quad\mbox{and}\quad j:\mathcal{Y}\lra B\Modcat$$ be the inclusions.
They are exact functors between exact categories and automatically induce derived functors between derived categories:
$$D(i):\D{\mathscr{E}}\lra \D{B}\quad\mbox{and}\quad D(j):\D{\mathcal{Y}}\lra \D{B}.$$
Moreover, by Lemma \ref{1-tilting}(1), there are another two additive functors between additive categories:
$$
\Hom_A(T, T\otimes_B-): \Pmodcat{B}\lra \mathcal{Y}:\;\; X\mapsto \Hom_A(T,T\otimes_BX)$$
$$\Coker(\eta'_{-}): \Pmodcat{B}\lra \mathscr{E}:\;\; X\mapsto \Coker(\eta'_X)
$$
for $X\in\Pmodcat{B}$. Here we fix an exact sequence
$$ 0\lra X \lraf{\eta_X'}\Hom_A(T, T\otimes_BX)\lra\Coker(\eta'_X)\lra 0 $$ for each $X$. Since $\K{B}_P$ is a triangulated subcategory of $\K{\Pmodcat{B}}$, we can define the following derived functors $\Psi: \D{B}\ra\D{\mathcal{Y}}$ and $\Phi: \D{B}\ra \D{\mathscr{E}}$ between derived categories of exact categories, where
$$
\xymatrix{\Psi:\;\; \D{B}\ar[r]^-{_p(-)}&\K{B}_P\ar[rr]^-{\Hom_A(T, T\otimes_B-)}&&\K{\mathcal{Y}}\ar[r]^-{Q_\mathcal{Y}}&\D{\mathcal{Y}}},
$$
$$
\xymatrix{\Phi:\;\; \D{B}\ar[r]^-{_p(-)}&\K{B}_P\ar[rr]^-{\Coker(\eta'_{-})}&&\K{\mathscr{E}}\ar[r]^-{Q_\mathscr{E}}&\D{\mathscr{E}}},
$$
with $Q_\mathcal{Y}$ and $Q_\mathscr{E}$ the localization functors.
By Lemma \ref{1-tilting}(2), there is a commutative diagram of natural transformations among triangle endofunctors of $\D{B}$:
$$
(\sharp)\quad
\xymatrix{
D(i)\circ\Phi [-1]\ar[r]\ar@{=}[d]&Id_{\D B}\ar[r]\ar@{=}[d]
&D(j)\circ\Psi\ar[d]^-{\simeq}\ar[r]& D(i)\circ\Phi\ar@{=}[d]\\
D(i)\circ\Phi [-1]\ar[r]&Id_{\D B}\ar[r]^-{\eta} & H\circ G\ar[r]&D(i)\circ\Phi}$$
This yields a commutative diagram of triangles in $\D{B}$ if applied to an object in $\D{B}$.

Now, let $D(H_T):\D{\mathcal{H}}\to \D{B}$ be the derived functor of $H_T$, and let $\mathbb{L}F:\D{B}\to \D{\mathcal{H}}$
be the total left-derived functor of $F$. Since $(F, H_T)$ is an adjoint pair by Lemma \ref{Baz}(3), we see that $(\mathbb{L}F, D(H_T))$ is an adjoint pair.

\begin{Lem}\label{FF}
$(1)$ Let $\overline{H_T}:\D{\mathcal H}\to \D{\mathcal Y}$ be the equivalence induced by $H_T$.
Then $$D(H_T)=D(j)\circ\overline{H_T},\;\; \Psi=\overline{H_T}\circ\mathbb{L}F\;\;\mbox{and}\;\; \Img(D(j))=\Img(H).$$
Moreover, $\big(\Psi, D(j)\big)$ is an adjoint pair and $D(j)$ is a fully faithful functor.

$(2)$ Let $\kappa:\Ker(G)\ra\D{B}$ and $\nu: \Img(H)\ra \D{B}$ be the inclusion functors. Then $D(i)=\kappa\circ \overline{D(i)}$ and $D(j)=\nu\circ \overline{D(j)}$, namely the commutative diagrams exist:
$$\xymatrix{\D{\mathscr{E}}\; \ar[rr]^{D(i)}\ar[dr]^{\overline{D(i)}}
                &  &    \D{B}    \\
                & \Ker(G)\quad\ar@{_{(}->}[ur]^{\kappa}                 } \; \; and \; \; \xymatrix{ \D{\mathscr{Y}}\; \ar[rr]^{D(j)}\ar[dr]^{\overline{D(j)}}
                &  &    \D{B}    \\
                & \Img(H)\; \; \ar@{_{(}->}[ur]^{\nu}                 }$$

$(3)$ $\Ker(D(i))=0$, $\Ker(\Phi)=\Img(H)$, and $\Phi$ commutes with direct products. Moreover, $\Phi$ induces a triangle functor
$$\overline{\Phi}:\;\;\D{B}/\Img(H)\lra\D{\mathscr{E}}$$ which commutes with direct products.

$(4)$ The two compositions
$$\Ker(G)\lraf{\kappa}\D{B}\lraf{Q}\D{B}/\Img(H)\quad\mbox{and}\quad
\D{B}/\Img(H)\lraf{\overline{\Phi}[-1]}\D{\mathscr{E}}\lraf{\overline{D(i)}}\Ker(G)$$
are quasi-inverse triangle equivalences, where $Q$ denotes the localization functor.
\end{Lem}

{\it Proof.} $(1)$  Clearly, $H_T$ induces an equivalence $H'_T: \mathcal{H}\to\mathcal{Y}$, its derived functor is denoted by $\overline{H}_T: \D{\mathcal{H}}\ra \D{\mathcal{Y}}$. Hence $H_T = j\circ H'_T$. Since all functors in this equality are exact, we have $D(H_T)=D(j)\circ\overline{H_T}$.
Note that $\mathbb{L}F$ is the composition of the functors:
$$
\xymatrix{\D{B}\ar[r]^-{_p(-)}&\K{B}_P\ar@{^{(}->}[r]&\K{\Pmodcat{B}}\,\ar@{^{(}->}[r]&\K{B}\ar[r]^-{F}&\K{\mathcal{H}}\ar[r]^-{Q_\mathcal{H}}&\D{\mathcal{H}}}.
$$
By Corollary \ref{Restriction}, the restriction of $F:\K{B}\to\K{\mathcal{H}}$ to $\K{\Pmodcat{B}}$ coincides with the one of ${_A}T\otimes_B-$, and has its image in
$\K{\Add({_A}T)}$, while the restriction of $H_T: \K{\mathcal{H}}\to\K{\mathcal{Y}}$ to $\K{\Add({_A}T)}$ coincides with the one of $\Hom_A(T,-)$.
It follows that $\Psi=\overline{H_T}\circ\mathbb{L}F$. Since $\overline{H_T}$ is an equivalence, $(\Psi, D(j))$ is an adjoint pair and $\Img(D(j))=\Img(D(H_T))$. By Lemma \ref{1-tilting}(3), $\Img(H)\subseteq \Img(D(H_T))$.

To show $\Img(D(H_T))\subseteq\Img(H)$, we apply the technique of homotopy limits in derived categories.

Let $\cpx{M}\in\C{B}$. By Lemma \ref{Homotopy}, $\cpx{M}\simeq \underleftarrow{\holim}(M^{(-n, n+1]})$ in $\D{B}$ where $n\in\mathbb{N}$. Suppose $\cpx{M}\in\C{\mathcal{Y}}$. Then $M^{(-n, n+1]}\in\Cb{\mathcal{Y}}$ by Corollary \ref{Restriction}(2). Further, by Theorem \ref{Torsion}, $\Img(H)$ is a triangulated subcategory of $\D{B}$ closed under direct products. It follows from $\mathcal{Y}\subseteq \Img(H)$ that $M^{(-n, n+1]}\in\Img(H)$ and further $\cpx{M}\in\Img(H)$.
This shows $\Img(D(H_T))\subseteq \Img(H)$, Thus $\Img(H)=\Img(D(H_T))=\Img(D(j))$.

Now, we show that the counit adjunction $\phi: \mathbb{L}F\circ D(H_T)\to \rm{Id}_{\D{\mathcal{H}}}$ is an isomorphism. This implies that both $D(H_T)$ and $D(j)$ are fully faithful functors
since $\overline{H_T}$ is an equivalence.

In fact, the adjoint pair $(F,H_T)$ of additive functors induces an adjoint pair $(K(F),K(H_T))$ of functors between $\K{\mathcal{H}}$ and $\K{B}$. Let $\phi': K(F)\circ K(H_T)\to Id_{\K{\mathcal H}}$ be its counit adjunction, and let $\psi': {\rm Id}_{\K{B}}\to K(H_T)\circ K(F)$ be its unit adjunction. Now, let $X\in\D{\mathcal{H}}$, $\cpx{Q}:={_p}(K(H_T)(X))$ and $\alpha_{\cpx{Q}}: \cpx{Q}\to K(H_T)(X)$ be the projective resolution of $K(H_T)(X)$ in $\D{B}$. Then $\phi_{X}:\mathbb{L}F\circ D(H_T)(X)\to X$ is given by the composition of $K(F)(\alpha_{\cpx{Q}})$ with $\phi'_{X}: K(F)\circ K(H_T)(X)\to X$. Since $H_T$ is fully faithful by Lemma \ref{Baz}(1), $K(H_T):\K{\mathcal{H}}\to \K{B}$ is fully faithful. This implies that $\phi'_X$ is an isomorphism by Lemma \ref{fact-adj}(1). Thus $\phi_{X}$ is an isomorphism in $\D{\mathcal{H}}$ if and only if $K(F)(\alpha_{\cpx{Q}})$ is a quasi-isomorphism. Since $H_T:\mathcal{H}\to B\Modcat$ is exact and fully faithful by Lemma \ref{Baz}(1), for each $U\in \mathscr{K}(\mathcal{H})$, we see that $U\in \mathscr{K}_{\rm {ac}}(\mathcal{H})$ if and only if $K(H_T)(U)\in\mathscr{K}_{\rm {ac}}(B)$. Consequently, $K(F)(\alpha_{\cpx{Q}})$ is a quasi-isomorphism if and only if so is $K(H_T)\circ K(F)(\alpha_{\cpx{Q}})$.
Since $\cpx{Q}\in\K{B}_P\subseteq\K{\Pmodcat{B}}$, we get $K(H_T) K(F)(\cpx{Q})=\Hom_A(T, T\otimes_B\cpx{Q})$ and $\psi'_{\cpx{Q}}=\eta'_{\cpx{Q}}$. Moreover,
$H\circ G(\cpx{Q})\simeq\Hom_A(T, T\otimes_B\cpx{Q})$ in $\D{B}$ by Lemma \ref{1-tilting}(2), and therefore there exists a diagram
$$
\xymatrix{
\cpx{Q}\ar@{=}[r]\ar[d]_-{\eta_{\cpx{Q}}}&\cpx{Q}\ar[rr]^-{\alpha_{\cpx{Q}}}\ar[d]_-{\psi_{\cpx{Q}}'} &&K(H_T)(X)\ar[d]^-{\simeq}_-{\psi_{H_T(X)}'}\\
H\circ G(\cpx{Q})\ar^-{\simeq}[r]&K(H_T)\circ K(F)(\cpx{Q})\ar[rr]^-{\;\;K(H_T)\circ K(F)(\alpha_{\cpx{Q}})\;\;} && K(H_T)\circ K(F)\circ K(H_T)(X).}$$
in which the left square is commutative in $\D{B}$ and the right square is commutative in $\K{B}$.
Note that the isomorphism $\psi_{H_T(X)}'$ follows from Lemma \ref{fact-adj}(2) and the fact that $K(H_T)(X)$ lies in $\Img(K(H_T))$. Since $\Img(H)=\Img(D(H_T))$ by $(1)$, we have $\cpx{Q}\simeq K(H_T)(X)\in\Img(H)$. Recall that $H$ is fully faithful.
Thus $\eta_{\cpx{Q}}$ is an isomorphism in $\D{B}$, and further, $K(H_T)\circ K(F)(\alpha_{\cpx{Q}})$ is an isomorphism in $\D{B}$.
As $K(H_T)\circ K(F)(\alpha_{\cpx{Q}})$ is represented by a chain map, it is a quasi-isomorphism. This shows that $\phi$ is a natural isomorphism.

$(2)$ This follows from $(1)$ and Lemma \ref{1-tilting}(3).

$(3)$ By Lemma \ref{thick}(2), $\cpx{X}\in\C{\mathscr{E}}$ is strictly exact if and only if it is exact in $\C{B}$.
This implies $\Ker(D(i))=0$, and therefore $\Ker(D(i)\circ \Phi)=\Ker(\Phi)$. By $(1)$ and $(\sharp)$,
we have $\Img(H)=\Ker(\Phi)$. Thus $\Phi$ induces a triangle functor $\overline{\Phi}: \D{B}/\Img(H)\ra \D{\mathscr{E}}$. Since $\Img(H)\subseteq\D{B}$ is closed under
direct products, $\D{B}/\Img(H)$ has direct products and the localization functor $Q$ in $(4)$ commutes with direct products
by the dual statement of \cite[Lemma 1.5]{BN}. As $\D{B}$ and $\D{B}/\Img(H)$ have the same objects,
$\overline{\Phi}$ commutes with direct products if and only if so does $\Phi$. Now let $\{M_i\}_{i\in I}$ be a family of
objects in $\D{B}$ together with the projections $\pi_i: \prod_{i\in I}M_i\to M_i$, where $I$ is a non-empty set.
To prove that the canonical morphism $\pi: \Phi(\prod_{i\in I}M_i)\to \prod_{i\in I}\Phi(M_i)$ in $\D{\mathscr{E}}$ induced from
$\{\Phi(\pi_i)\}_{i\in I}$ is an isomorphism, it suffices to show that $D(i)(\pi)$ is an isomorphism in $\D{B}$. The reason reads as follows: $\pi$ can be embedded into a triangle $\Phi(\prod_{i\in I}M_i)\to \prod_{i\in I}\Phi(M_i)\ra W\ra \Phi(\prod_{i\in I}M_i)[1]$ in $\D{\mathscr{E}}$, and if $D(i)(\pi)$ is an isomorphism, then $D(i)(W)=0$. It follows from $\Ker(D(i))=0$ that $W=0$ and $\pi$ is an isomorphism in $\D{\mathscr{E}}$.

Since $D(i)$ commutes with direct products by Lemma \ref{thick}(1), it is enough to show that
the canonical morphism $D(i)\circ \Phi(\prod_{i\in I}M_i)\to \prod_{i\in I}(D(i)\circ \Phi(M_i))$ in $\D{B}$
induced from  $\{D(i)\circ \Phi(\pi_i)\}_{i\in I}$ is an isomorphism, that is, $D(i)\circ\Phi$ commutes with direct products.
But this follows from the diagram $(\sharp)$ and the fact that $H$ and $G$ commute with direct products. Thus $\Phi$ and therefore $\overline{\Phi}$ commute with direct products.

$(4)$ By Theorem \ref{Torsion}, the pair $(\Ker(G),\Img(H))$ is a semi-orthogonal decomposition of $\D{B}$. By $(3)$, $\Ker(\overline{D(i)}\circ \Phi[-1])=\Ker(\Phi)=\Img(H)$. It follows from
the diagram $(\sharp)$ and Lemma \ref{Torsion pair} that $\overline{D(i)}\circ \Phi[-1]: \D{B}\to \Ker(G)$ is a right adjoint of $\kappa$
and induces a triangle equivalence $\overline{D(i)}\circ \overline{\Phi}[-1]:\D{B}/\Img(H)\to\Ker(G)$. Thus $(4)$ holds. $\square$
\medskip

As a consequence of Theorem \ref{Torsion} and Lemma \ref{FF}(1), we re-obtain \cite[Theorem 4.5]{Bz1} which was proved by employing the theory of derivators in \cite{STK}.

\begin{Koro}
There is a triangle equivalence $\D{\mathcal H}\to \D{A}$ given by the composition of the equivalence $D(H_T): \D{\mathcal H}\to\Img(H)$ with
a quasi-inverse $\Img(H)\to \D{A}$ of the equivalence $H: \D{A}\to \Img(H)$.
\end{Koro}

\begin{Rem}\label{Bounded-above}
The functors $H$ and $G$ can be restricted to functors between bounded-above derived categories since both ${_A}T$ and $T_B$ have finite projective dimension.
This implies that $\big(\Df{B}\cap\Ker(G), \Df{B}\cap\Img(H)\big)$ is a semi-orthogonal decomposition of $\Df{B}$.
Moreover, by definition, $\Phi$ and $\Psi$ can also be restricted to functors between bounded-above derived categories:
$$\Phi^-:\;\;\Df{B}\to \Df{\mathscr{E}}\quad \mbox{and}\quad \Psi^-:\;\;\Df{B}\to \Df{\mathcal Y}.$$
Thus Lemma \ref{FF} hold true for bounded-above derived categories.
\end{Rem}

\subsection{Fully faithful triangle functors between derived categories\label{subsect3.3}}

In this section we develop properties of triangle functors needed in our proofs.

By the diagram $(\sharp)$, there exists a triangle in $\D{B}$ for each $M\in B\Modcat$:
$$
(\diamondsuit)\quad
D(i)\circ\Phi [-1](M)\lra M\lra D(j)\circ\Psi(M)\lra  D(i)\circ\Phi(M),
$$
with $\Psi(M)\in \Cf{\mathcal{Y}}$, $\Phi(M)\in\Cf{\mathscr{E}}$ and $\Psi(M)^i=0=\Phi(M)^i$ for all $i\geq 1$.
Taking homologies of the triangle, we get a $5$-term exact sequence of $B$-modules:
$$
\varepsilon_M:\quad 0\lra Y_M\lraf{\epsilon^{-2}_M} X_M\lraf{\epsilon^{-1}_M} M\lraf{\epsilon^{0}_M} Y^M\lra X^M\lra 0,
$$
where $ Y_M:=H^{-1}\big(D(j)\circ\Psi(M)\big),\; X_M:=H^{-1}\big(D(i)\circ\Phi(M)\big),
\; Y^M:=H^0\big(D(j)\circ\Psi(M)\big)$ and $X^M:=H^0\big(D(i)\circ\Phi(M)\big).$ In the following, we consider the full subcategories of $B\Modcat$:
$$\mathcal{E}:=\{M\in B\Modcat\mid X^M=0\},$$
$$\mathcal{X}:=\{X\in B\Modcat\mid \Ext_B^n(X,Y)=0,\, Y\in\mathcal{Y},\, n=0,1\},$$
$${^\bot}\mathcal{Y}:=\{X\in B\Modcat\mid \Ext_B^n(X,Y)=0,\, Y\in\mathcal{Y},\, n\geq 0 \},$$
$$\mathscr{E}^\bot:=\{Y\in B\Modcat\mid\Ext_B^n(X, Y)=0,\, X\in\mathscr{E},\, n\geq 0\},$$
$$\Coker(\mathscr{E}):=\{\Coker(f)\in B\Modcat\mid f:E_1\to E_0\;\mbox{with}\; E_1, E_0\in\mathscr{E}\}.$$
Then $\mathcal{X}$
is closed under extensions and cokernels in $B\Modcat$, and both $\mathscr{E}^\bot$ and ${^\bot}\mathcal{Y}$ are thick subcategories of $B\Modcat$. Moreover, we have the following result which generalizes \cite[Lemma 5.8]{Bz1}.

\begin{Lem}\label{Ext-orthogonal}
$\mathscr{E}^\bot=\mathcal{Y}$ and $\mathscr{E}={^\bot}\mathcal{Y}$.
\end{Lem}

{\it Proof.} We have $\mathscr{E}\subseteq \Ker(G)$ and $\mathcal{Y}\subseteq \Img(H)$.
Since $\Hom_{\D{B}}(\cpx{X}, \cpx{Y})=0$ for any $\cpx{X}\in\Ker(G)$ and $\cpx{Y}\in\Img(H)$ by Theorem \ref{Torsion},
we have $\mathcal{Y}\subseteq\mathscr{E}^\bot$ and $\mathscr{E}\subseteq{^\bot}\mathcal{Y}$.
To show $\mathscr{E}^\bot\subseteq\mathcal{Y}$, we first show that if $Z\in \mathscr{E}^\bot$, then $\Hom_{\D{B}}(\cpx{E}, Z)=0$ for all $\cpx{E}\in\Cf{\mathscr{E}}$.

Indeed, let $\mathscr{D}_Z$ be the full subcategory of $\D{B}$ consisting of all complexes
$\cpx{U}$ such that $\Hom_{\D B}(\cpx{U}[n], Z)$ = $0$ for all $n\in\mathbb{Z}$. Then $\mathscr{D}_Z$ is a triangulated subcategory
of $\D{B}$ which is closed under direct sums and contains $\mathscr{E}$. Thus homotopy colimits in $\D B$ of sequences in $\mathscr{D}_Z$ belong to $\mathscr{D}_Z$ by Definition \ref{def-holim}. By the first isomorphism in Lemma \ref{Homotopy},
each complex in $\Cf{\mathscr{E}}$ can be obtained in $\D{B}$ from bounded complexes in $\Cb{\mathscr{E}}$
by taking homotopy colimits. Hence $\Cf{\mathscr{E}}\subseteq \mathscr{D}_Z$, and therefore $\Hom_{\D{B}}(\cpx{E}, Z)=0$ for all $\cpx{E}\in\Cf{\mathscr{E}}$.

Since $\Phi(Z)[-1]\in\Cf{\mathscr{E}}$, we have $\Hom_{\D B}(\Phi(Z)[-1], Z)=0$. It follows from the triangle
$(\diamondsuit)$ that $D(j)\circ\Psi(Z)\simeq Z\oplus D(i)\circ\Phi(Z)$.
As $\Hom_{\D B}(D(i)\circ\Phi(Z), D(j)\circ\Psi(Z))=0$, we get $D(i)\circ\Phi(Z)=0$ and $Z\simeq D(j)\circ\Psi(Z)\in\Img(H)$.
Thus $Z\in\Img(H)\cap B\Modcat=\mathcal{Y}$ by Lemma \ref{Baz}(2). This shows $\mathscr{E}^\bot\subseteq\mathcal{Y}$, and therefore $\mathscr{E}^\bot=\mathcal{Y}$.

Since $\mathcal{Y}$ is an abelian category of $B\Modcat$ by Corollary \ref{Restriction}(2), it is closed under cokernels in $B\Modcat$.
By the last isomorphism in Lemma \ref{Homotopy}, each complex in $\Cf{\mathscr{Y}}$ can be obtained
in $\D{B}$ from bounded complexes in $\Cb{\mathscr{Y}}$ by taking homotopy limits.
Similarly, we show that if $Z\in {^\bot}\mathcal{Y}$, then $\Hom_{\D{B}}(Z, \cpx{Y})=0$ for any $\cpx{Y}\in\Cf{\mathcal{Y}}$.
Then the inclusion ${^\bot}\mathcal{Y}\subseteq \mathscr{E}$ follows from $(\diamondsuit)$ and $\Psi(Z)\in \Cf{\mathcal{Y}}$.
Thus $\mathscr{E}={^\bot}\mathcal{Y}$. $\square$

\medskip
Next, we study properties of $\varepsilon_M$ and $\mathcal{E}$.

\begin{Lem} \label{B-cat}
The following statements hold true for $M\in B\Modcat$.

$(1)$ $Y_M, Y^M\in\mathcal{Y}$ and $X^M\in\Coker(\mathscr{E})$. If $M\in\mathcal{E}$, then $X_M\in\Coker(\mathscr{E})$.

$(2)$ $\mathcal{X}=\Coker(\mathscr{E})=\{M\in\mathcal{E}\mid Y_M=0=Y^M\}$
and $\mathcal{Y}=\{M\in\mathcal{E}\mid X_M=0\}$.

$(3)$ $\mathcal{E}$ is closed under extensions and quotients in $B\Modcat$.
\end{Lem}

{\it Proof.} $(1)$ Since $\mathcal{Y}$ is an abelian category of $B\Modcat$ by Corollary \ref{Restriction}(2), it is clear that $H^i(\Psi(M))\in\mathcal{Y}$ for any $i\in\mathbb{Z}$. Moreover, there is a right exact sequence $\Phi(M)^{-1}\lraf{d^{-1}} \Phi(M)^0\to X^M\to 0$ with $\Phi(M)^{-1}, \Phi(M)^0\in\mathscr{E}$.
This shows $X^M\in\Coker(\mathscr{E})$.

Let $M\in \mathcal{E}$. Then $X^M=0$, which means that the homomorphism $d^{-1}$ is surjective.
Since $\mathscr{E}$ is closed under kernels of surjective homomorphisms in $B\Modcat$ by Lemma \ref{thick}(1), we have
$\Ker(d^{-1})\in\mathscr{E}$. Similarly, there is a right exact sequence $\Phi(M)^{-2}\to \Ker(d^{-1})\to X_M\to 0$ such that $\Phi(M)^{-2}\in\mathscr{E}$.
Thus $X_M\in\Coker(\mathscr{E})$.

$(2)$ Since $\mathcal{X}$ is closed under cokernels in $B\Modcat$ and $\mathscr{E}\subseteq \mathcal{X}$, we have $\Coker(\mathscr{E})\subseteq \mathcal{X}$.

Let $M\in\mathcal{X}$. Since $Y^M\in\mathcal{Y}$ by $(1)$, $\Hom_B(M, Y^M)=0$ and $Y^M\simeq X^M$.
As $X^M$ is a quotient of $\Phi(M)^0\in\mathscr{E}$, it follows from Lemma \ref{Ext-orthogonal} that $\Hom_B(X^M, Y^M)=0$.
This forces $X^M=0=Y^M$, and thus $M\in \mathcal{E}$. Since $Y_M\in\mathcal{Y}$, we have
$\Ext_B^1(M, Y_M)=0$. Then $M\simeq Y_M\oplus X_M$. Note that $\Hom_B(M, Y_M)=0$ for $M\in\mathcal{X}$ and $Y_M\in\mathcal{Y}$. Consequently, $Y_M=0$ and $M\simeq X_M\in\Coker(\mathscr{E})$. Thus $\mathcal{X}\subseteq\{M\in\mathcal{E}\mid Y_M=0=Y^M\}$.
Conversely, if $N\in\mathcal{E}$ and $Y_N=0=Y^N$, then $N\simeq X_N\in\Coker(\mathscr{E})$ by $(1)$.

If $N\in\mathcal{E}$ with $X_N=0$, then $N\simeq Y^N\in\mathcal{Y}$. Since $\mathcal{Y}\subseteq \Img(H)$, it follows from Lemma \ref{FF}(1) that $N\in\mathcal{Y}$ if and only if $\Phi(N)=0$ in $\D{B}$. This implies the last equality in $(2)$.

$(3)$ Let $0\to U\to V\to W\to 0$ be an exact sequence of $B$-modules. Then this sequence induces a triangle
$$D(i)\circ\Phi(U)\lra D(i)\circ\Phi(V)\lra D(i)\circ\Phi(W)\lra D(i)\circ\Phi(U)[1]$$ in $\D{B}$. Taking homology $H^0$ of this triangle yields an exact sequence
$X^U\to X^V\to X^W\to H^1(D(i)\circ \Phi(U))$ in $B\Modcat$. By the definition of $\Phi$, $\Phi(U)$ is a complex with $0$ at all positive degree, this implies $H^1(D(i)\circ \Phi(U))=0$. Hence $\mathcal{E}$ is closed under
extensions and quotients in $B\Modcat$. $\square$

\begin{Koro}\label{adjoint}
$(1)$ For each $M\in\mathcal{E}$, there exists a $4$-term exact sequence of $B$-modules:
$$ 0\lra Y_M\lraf{\varepsilon_M^{-2}} X_M\lraf{\varepsilon_M^{-1}} M
\lraf{\varepsilon_M^{0}} Y^M\lra 0$$ with $Y_M, Y^M\in\mathcal{Y}$ and $X_M\in\mathcal{X}$.

$(2)$ $\varepsilon_M^{-1}$ and $\varepsilon_M^{0}$ in (1) give rise to isomorphisms of abelian groups for $X\in\mathcal{X}$ and $Y\in\mathcal{Y}$:
$$
(\varepsilon_M^{-1})^*: \Hom_B(X, X_M)\lraf{\simeq} \Hom_B(X, M)\quad\mbox{and}\quad
(\varepsilon_M^{0})_*: \Hom_B(Y^M, Y)\lraf{\simeq} \Hom_B(M, Y).
$$

$(3)$ The inclusion $\mathcal{X}\to \mathcal{E}$ has a right adjoint $r:\mathcal{E}\to \mathcal{X}$
given by $M\mapsto X_M$ for any $M\in\mathcal{E}$. Moreover, if $0\to M_1\to M_2\to M_3\to 0$
is an exact sequence in $B\Modcat$ with $M_i\in\mathcal{E}$ for $1\leq i\leq 3$, then $r(M_1)\to r(M_2)\to r(M_3)\to 0$ is exact.
\end{Koro}

{\it Proof.} $(1)$ We have an exact sequence $\epsilon_M$. Since $M\in\mathcal{E}$, $\epsilon_M$ can be shortened into $4$-term sequence with the desired properties by Lemma \ref{B-cat}(1)-(2).

(2) Since $X\in \mathcal{X}$ and $Y_M, Y^M\in\mathcal{Y}$,
we have $\Hom_B(X, Y^M)=\Hom_B(X, Y_M)=\Ext_B^1(X, Y_M)=0$. It then follows from the exact sequence in (1) that $(\varepsilon_M^{-1})^*$ is an isomorphism. This
implies the first isomorphism in $(2)$. Similarly, $(\varepsilon_M^{0})_*$ is an isomorphism.

(3) If $M\in\mathcal{X}$, then
$\varepsilon_M^{-1}$ is an isomorphism by the proof of Lemma \ref{B-cat}(2).
Now, the isomorphism $(\varepsilon_M^{-1})^*$ in $(2)$ gives an adjunction isomorphism of adjoint pairs of the functors.

Applying the functor $D(i)\circ\Phi:\D{B}\to\D{B}$ to the exact sequence $0\to M_1\to M_2\to M_3\to 0$ in $B\Modcat$, we obtain a triangle
$$D(i)\circ\Phi(M_1)\lra D(i)\circ\Phi(M_2)\lra D(i)\circ\Phi(M_3)\lra D(i)\circ\Phi(M_1)[1].$$
If we take $H^{-1}$ on the triangle, then we gain an exact sequence $r(M_1)\to r(M_2)\to r(M_3)\to X^{M_1}$. Now, it follows from $M_1\in\mathcal{E}$ that $X^{M_1}=0$, and therefore the second part of $(3)$ holds. $\square$
\medskip

The following property of $\mathcal{E}$ is crucial in establishing an equivalence of derived categories of exact categories in Lemma \ref{FFEX} below.

\begin{Lem}\label{InjectiveinE}
If $M$ is an injective $B$-module, then $M\in\mathcal{E}$.
\end{Lem}

To show Lemma \ref{InjectiveinE}, we need to determine the image of an injective cogenerator for $B\Modcat$ under $H\circ G$.
Let $$(-)^{\vee}:=\Hom_{\mathbb{Z}}(-, \mathbb{Q}/\mathbb{Z}):\;\; \mathbb{Z}\Modcat\lra\mathbb{Z}\Modcat.$$
Then $(-)^\vee$ induces an exact functor $B\opp\Modcat\to B\Modcat$ and $B^\vee$ is an injective cogenerator for $B\Modcat$.
Further, we define two natural transformations:
$$
\theta: T\otimes_B\Hom_A(-,T)^\vee\to \Hom_A(-, A)^\vee:\;\;  A\Modcat\lra A\Modcat,
$$
$$
\rho: \Hom_A(-,A)\otimes_AT\lra \Hom_A(-, T):\;\; A\Modcat\to B\opp\Modcat
$$
given by
$$\theta_X:\; {_A}T\otimes_B\Hom_A(X,T)^\vee\lra\Hom_A(X, A)^\vee:\;\; t\otimes\sigma\mapsto [f\mapsto (f (\cdot t))\sigma],$$
$$\rho_X: \Hom_A(X,A)\otimes_AT\lra \Hom_A(X, T):\;\;f\otimes t\mapsto f (\cdot t),$$
where $X\in A\Modcat, t\in T, \sigma\in \Hom_A(X,T)^\vee, f\in \Hom_A(X,A)$ and $(\cdot t)\in\Hom_A(A,T)$
is the right multiplication by $t$.

\begin{Lem} \label{Basic}
$(1)$ $\theta$ is a natural isomorphism.

$(2)$ If $X\in\Pmodcat{A}$, then

$\quad (i)$ $\Tor_j^B(T, \Hom_A(X,T)^\vee)=0$ for all $j\geq 1$.

$\quad (ii)$ $\rho_X$ is injective and
$\Tor_j^A(\Hom_A(X, A), T)=0$ for all $j\geq 1$.
\end{Lem}

{\it Proof.}
$(1)$ For $U\in A\Modcat$ and $g\in\Hom_A(U, T)$, we define
$$\theta_{U,\, X}:\; \Hom_A(U, T)\otimes_B\Hom_A(X,T)^\vee\lra  \Hom_A(X, U)^\vee:\;\; g\otimes \sigma\mapsto [f\mapsto (fg)\sigma]$$
Clearly, $\theta_{A,\, X}=\theta_X$ under the identification of $\Hom_A(A,T)$ with $T$.
Moreover, if $U\in\add(T)$, then $\theta_{U,\, X}$ is an isomorphism. Note that $T_B$ has a finitely generated projective
resolution of length at most $n$:
$$
(\ast)\quad 0\lra \Hom_A(T_n, T)\lra \cdots \lra\Hom_A(T_1, T)\lra\Hom_A(T_0,
T)\lra T_B\lra 0$$ with  $T_j\in\add(_AT)$ for all $0\leq j\leq n$.
Now, we can construct the following exact commutative diagram:
$$
\xymatrix{
\Hom_A(T_1, T)\otimes_B\Hom_A(X,T)^\vee\ar[r]\ar^-{\theta_{T_1,\, X}}[d]_-{\simeq}&\Hom_A(T_0, T)\otimes_B\Hom_A(X,T)^\vee\ar^-{\theta_{T_0,\, X}}[d]_-{\simeq}\ar[r]& T\otimes_B\Hom_A(X,T)^\vee\ar^-{\theta_X}[d]\ar[r]&0\\
\Hom_A(X, T_1)^\vee\ar[r]& \Hom_A(X, T_0)^\vee\ar[r]&\Hom_A(X,A)^\vee\ar[r]&0.}$$
Thus $\theta_X$ is an isomorphism. This shows $(1)$.

(2) If $X$ is projective, then the sequence $$0\lra \Hom_A(X, T_n)^\vee\lra\cdots\lra \Hom_A(X, T_1)^\vee\lra\Hom_A(X, T_0)^\vee\lra\Hom_A(X,A)^\vee\lra 0$$ is exact.
Now $(i)$ follows if we apply $-\otimes_B\Hom_A(X,T)^\vee$ to $(\ast)$.

Note that $(ii)$ holds if and only if the canonical map $A^{\alpha}\otimes_AT\to T^\alpha$ is injective and $\Tor^A_j(A^\alpha, T)=0$ for any nonempty set $\alpha$ and for any $j\geq 1$. The latter follows from \cite[Lemmas 2.5 and 2.4(3)]{XC2}. $\square$

\begin{Lem}\label{Injective}
Let $\cpx{P}: 0\to P_n\to \cdots\to P_1\to P_0\to 0$ be the deleted projective resolution of ${_A}T$ (see Definition \ref{Tilting}), where $P_i$ are of degree $-i$ for $0\leq i\leq n$. Then

$(1)$ $G(B^\vee)\simeq \Hom_A(\cpx{P}, A)^\vee$ in $\D{A}$.

$(2)$ Let $\cpx{\rho}:=(\rho_{P_i})_{0\leq i\leq n}:\;\;\Hom_A(\cpx{P}, A)\otimes_AT\lra \Hom_A(\cpx{P}, T).$
Then $\cpx{\rho}$ is an injective chain map and there is a commutative diagram of triangles in $\D{B}$:
$$
\xymatrix{
\Coker(\cpx{\rho})^\vee\ar[r]\ar@{=}[d]&\Hom_A(\cpx{P}, T)^\vee\ar[r]^-{\cpx{\rho}{^\vee}}\ar[d]^-{\simeq}
&(\Hom_A(\cpx{P}, A)\otimes_AT)^\vee\ar^-{\simeq}[d]\ar[r]& \Coker(\cpx{\rho})^\vee[1]\ar@{=}[d]\\
\Coker(\cpx{\rho})^\vee\ar[r]&B^\vee\ar[r]^-{\eta_{B^\vee}} & H\circ G(B^\vee)\ar[r]&\;\; \; \Coker(\cpx{\rho})^\vee[1].}$$
\end{Lem}

{\it Proof.} $(1)$ By $(T1)$ and $(T2)$ in Definition \ref{Tilting}, there is an exact sequence of $B$-modules
$$
0\lra \Hom_A(P_n, T)^\vee\lra \cdots \lra\Hom_A(P_1, T)^\vee\lra\Hom_A(P_0,T)^\vee\lra B^\vee\lra 0,$$
where $P_i\in \Pmodcat{A}$ for $0\leq i\leq n$. This means that $\Hom_A(\cpx{P}, T)^\vee$ is quasi-isomorphic to $B^\vee$. Therefore $G(B^\vee)\simeq G(\Hom_A(\cpx{P}, T)^\vee)$ in $\D{A}$. Moreover, by Lemma \ref{Basic}(2), $\Tor_j^B(T, \Hom_A(P_i,T)^\vee)=0$ for any $j\geq 1$. This yields
$G(\Hom_A(\cpx{P}, T)^\vee)\simeq T\otimes_B \Hom_A(\cpx{P}, T)^\vee$ by Lemma \ref{HT}(1). Since
$T\otimes_B \Hom_A(\cpx{P}, T)^\vee\simeq \Hom_A(\cpx{P}, A)^\vee$ by Lemma \ref{Basic}(1), we have $G(B^\vee)\simeq \Hom_A(\cpx{P}, A)^\vee$ in $\D{A}$.

$(2)$ For $X\in\Pmodcat{A}$ and $j\geq 1$, we have $\Ext_A^j(T, \Hom_A(X,A)^\vee)\simeq \Tor_j^A(\Hom_A(X, A), T)^\vee=0$ by Lemma \ref{Basic}(3).
Since $P_i\in\Pmodcat{A}$ for all $0\leq i\leq n$, Lemma \ref{HT}(2) implies that
$$\Hom_A(T, \beta_{\Hom_A(\cpx{P}, A)^\vee}):\;\; \Hom_A(T, \Hom_A(\cpx{P}, A)^\vee)\lra H(\Hom_A(\cpx{P}, A)^\vee)$$
is an isomorphism in $\D{B}$. Consequently, there are isomorphisms in $\D{B}$:
$$H\circ G(B^\vee)\simeq H(\Hom_A(\cpx{P}, A)^\vee)\simeq \Hom_A(T, \Hom_A(\cpx{P}, A)^\vee)\simeq (\Hom_A(\cpx{P}, A)\otimes_AT)^\vee.$$

Clearly, $\cpx{\rho}$ is injective by Lemma \ref{Basic}(3). Thus the sequence $$0\lra \Hom_A(\cpx{P}, A)\otimes_AT\lraf{\cpx{\rho}} \Hom_A(\cpx{P}, T)\lra \Coker(\cpx{\rho})\lra 0$$ is exact, and yields an triangle in $\D{B^{op}}$
$$ \Coker(\cpx{\rho})[-1]\lra \Hom_A(\cpx{P}, A)\otimes_AT\lraf{\cpx{\rho}} \Hom_A(\cpx{P}, T)\lra \Coker(\cpx{\rho})$$ and a triangle in $\D{B}$
$$ \Coker(\cpx{\rho})^\vee\lra \Hom_A(\cpx{P}, T)^\vee\lraf{\cpx{\rho}{^\vee}} (\Hom_A(\cpx{P}, A)\otimes_AT)^\vee\lra \Coker(\cpx{\rho})^\vee[1].$$
Now, it follows from this triangle that (2) can be deduced by the commutative diagram in $\D{B}$:
$$
\xymatrix{
\Hom_A(\cpx{P}, T)^\vee\ar[rr]^-{\eta'_{\Hom_A(\cpx{P}, T)^\vee}}\ar[d]^-{\simeq}
&&\Hom_A\big(T, T\otimes_B \Hom_A(\cpx{P}, T)^\vee\big)\ar^-{\simeq}[d]\ar[r]^-{\simeq}&(\Hom_A(\cpx{P}, A)\otimes_AT)^\vee\ar[d]\\
B^\vee\ar[rr]^-{\eta_{B^\vee}} && H\circ G(B^\vee)\ar@{=}[r]&\;\; \; H\circ G(B^\vee)}$$
where the composition of the first arrow is $\cpx{\rho}{^\vee}$. $\square$

\medskip
{\bf Proof of Lemma \ref{InjectiveinE}}: Let $M$ be an injective $B$-module. Then $M\in\Prod(B^\vee)$. Since  $X^M = H^0(D(i)\circ \Phi(M))$ and since
the functors $D(i)$, $\Phi$ and $H^0$ commute with direct products, we have $X^M\in\Prod(X^{B^\vee})$. By the diagram $(\sharp)$ and Lemma \ref{Injective}(2), $D(i)\circ\Phi(B^\vee)\simeq \Coker(\cpx{\rho})^\vee[1]$ in $\D{B}$. Since $H^0\big(\Coker(\cpx{\rho})^\vee[1]\big)=\big(H^{-1}(\Coker(\cpx{\rho}))\big)^\vee=0$, we have  $X^{B^\vee}=0$. Thus $X^M=0$. $\square$

\begin{Lem}\label{FFEX}
$(1)$ The inclusion $\mathcal{E}\subseteq B\Modcat$ induces a triangle equivalence $\D{\mathcal{E}}\ra \D{B}$ which restricts to an equivalence $\mathscr{D}^*(\mathcal{E})\ra \mathscr{D}^*(B)$ for any $*\in\{+, -, b\}$.

$(2)$ The inclusion $\mathscr{E}\subseteq \mathcal{E}$ induces a fully faithful triangle functor $\Df{\mathscr{E}}\to\Df{\mathcal{E}}$.
\end{Lem}

{\it Proof.}
$(1)$ By Lemma \ref{B-cat}(3), $\mathcal{E}$ is closed under extensions and quotients in $B\Modcat$. In particular,
$\mathcal{E}$ is a fully exact subcategory of $B\Modcat$, which is closed under direct summands. By Lemma \ref{Injective}, $\mathcal{E}$ contains all injective $B$-modules.
Thus $\mathcal{E}$ satisfies the two properties:

$(i)$ If $0\to X\to Y\to Z\to 0$ is exact in $B\Modcat$ with $Y\in\mathcal{E}$, then $Z\in\mathcal{E}$.

$(ii)$ For each $B$-module $M$, there exists a short exact sequence $0\to M\to E_0\to E_1\to  0$ in $B\Modcat$ such that $E_0$ is injective
and $E_1\in\mathcal{E}$.

Now, $(1)$ follows from the dual statement of Lemma \ref{resolution}(2) with $\mathcal{A}=B\Modcat$.

$(2)$ By Lemma \ref{thick}(1), if $0\to X\to Y\to Z\to 0$ is exact in $B\Modcat$ with $Z\in \mathscr{E}$, then $X\in\mathscr{E}$ if and only if $Y\in\mathscr{E}$. This implies that $\mathscr{E}$ is a fully exact subcategory of $B\Modcat$. Clearly,
$\mathscr{E}\subseteq\mathcal{X}\subseteq \mathcal{E}$ by Lemma \ref{Ext-orthogonal} and Lemma \ref{B-cat}(2).

Let $0\to E_1\to E_0\lraf{g} F\to 0$ be an exact sequence in $B\Modcat$ with $E_1, E_0\in\mathcal{E}$ and $F\in\mathscr{E}$.
By Corollary \ref{adjoint}(3), we can apply the functor $r:\mathcal{E}\to\mathcal{X}$ to $g$ and obtain a surjective map $r(g):r(E_0)\to r(F)$.
Moreover, the composition of $r(g)$ with the counit map $\varepsilon_F^{-1}: r(F)\to F$ coincides with
the composition of the counit map $\varepsilon_{E_0}^{-1}:r(E_0)\to E_0$ with $g$.
Note that $\varepsilon_F^{-1}$ is an isomorphism for $F\in \mathscr{E}\subseteq\mathcal{X}$.
Since $r(E_0)\in\mathcal{X}=\Coker(\mathscr{E})$ by Lemma \ref{B-cat}(2), there exists a surjective map $f:F_0\to r(E_0)$ with $F_0\in\mathscr{E}$.
Let $h:=f\varepsilon_{E_0}^{-1}$. Then $hg=f\varepsilon_{E_0}^{-1}g=fr(g)\varepsilon_F^{-1}:F_0\to F$, and therefore $hg$ is surjective.
Let $F_1:=\Ker(hg)$. Since $F_0, F\in\mathscr{E}$, we have $F_1\in\mathscr{E}$.
Moreover, there is an exact commutative diagram:
$$
\xymatrix{0\ar[r]&F_1\ar[r]\ar[d]&F_0\ar[r]^-{hg}\ar[d]_-{h}& F\ar[r]\ar@{=}[d]&0\\
0\ar[r]&E_1\ar[r]& E_0\ar[r]^-{g}& F\ar[r]&0.
}
$$
Now, $(2)$ follows from Lemma \ref{resolution}(1). $\square$

\medskip
Finally, we prove the following  result which is crucial to the proof of Theorem \ref{Main result}.
\begin{Prop}\label{Embedding}
The functor $\overline{D(i)}:\D{\mathscr{E}}\to \Ker(G)$ is a triangle equivalence.
\end{Prop}

{\it Proof.} Since $D^{-}(i):\Df{\mathscr{E}}\to \Df{B}$ is the composition of the induced functor $\Df{\mathscr{E}}\to\Df{\mathcal{E}}$ with the one $\Df{\mathcal{E}}\to\Df{B}$, Lemma \ref{FFEX} implies that $D^{-}(i)$ is fully faithful. By Lemmas \ref{1-tilting}(3) and \ref{thick}(2), $\Img(D^{-}(i))=\Df{B}\cap \Ker(G)$.
Consequently, the restriction
$$\overline{D^-(i)}:\Df{\mathscr{E}}\lra \Df{B}\cap \Ker(G)$$
of the functor $\overline{D(i)}$ is a triangle equivalence. Moreover, by Lemma \ref{FF}(4) and Remark \ref{Bounded-above},
the two compositions:
$$(\ast\ast) \quad \Df{B}\cap\Ker(G)\lraf{\kappa^-}\Df{B}\lraf{Q^-}\Df{B}/(\Df{B}\cap\Img(H))\quad\mbox{and}$$
$$
\Df{B}/(\Df{B}\cap\Img(H))\lraf{\overline{\Phi^-}[-1]}\Df{\mathscr{E}}\lraf{\overline{D^-(i)}}\Df{B}\cap\Ker(G)$$
are quasi-inverse triangle equivalences, where $^-$ denotes universally the restriction of the involved functors to bounded-above derived categories. This implies that $\overline{\Phi^-}[-1]$ is a triangle equivalence.

Next, we show that the inclusion $\Df{B}\to\D{B}$ induces a fully faithful functor
$$\Df{B}/(\Df{B}\cap\Img(H))\lra \D{B}/\Img(H).$$
By the dual statement of \cite[Lemma 10.3]{Keller}, it suffices to show that each morphism $\cpx{f}: \cpx{M}\to \cpx{Y}$
in $\D{B}$ with $\cpx{M}\in\Df{B}$ and $\cpx{Y}\in\Img(H)$ factorizes through an object $\cpx{Z}\in\Df{B}\cap\Img(H)$.
Since $\Kf{\Pmodcat{B}}$ is equivalent to $\Df{B}$, we can assume that $\cpx{M}\in\Cf{\Pmodcat{B}}$ and $\cpx{f}$ is represented
by a chain map. As $\Img(H)=\Img(D(j))$ by Lemma \ref{FF}(1), we have $\cpx{Y}\in\C{\mathcal{Y}}$ up to isomorphism in $\D{B}$.
Then $\cpx{f}$ factorizes through the left truncated complex $\tau^{\leq m}\cpx{Y}:\; \cdots \to Y^i\to Y^{i+1}\to\cdots \to Y^{m-1}\to \Ker(d_Y^{m})\to 0$ of $\cpx{Y}$ at some degree $m$. Since $\mathcal{Y}$ is an abelian subcategory of $B\Modcat$, $\Ker(d_Y^{m})\in\mathcal{Y}$ and $\tau^{\leq m}\cpx{Y}\in\Cf{\mathcal{Y}}$, and therefore $\tau^{\leq m}\cpx{Y}\in\Df{B}\cap\Img(H)$.

Now, we construct the commutative diagram of functors:
$$
\xymatrix{
\Df{\mathscr{E}}\ar[r]^-{\overline{D^-(i)}}_-{\simeq}\ar@{_{(}->}[d]&\Df{B}\cap\Ker(G)\ar[r]^-{\kappa^-}\ar@{_{(}->}[d]&\Df{B}\ar[r]^-{Q^-}\ar@{_{(}->}[d]
&\Df{B}/(\Df{B}\cap\Img(H))\ar[r]^-{\overline{\Phi^-}[-1]}_-{\simeq}\ar@{_{(}->}[d]&\Df{\mathscr{E}}\ar@{_{(}->}[d]\\
\D{\mathscr{E}}\ar[r]^-{\overline{D(i)}}&\Ker(G)\ar[r]^-{\kappa}&\D{B}\ar[r]^-{Q}&\D{B}/\Img(H)
\ar[r]^-{\overline{\Phi}[-1]}&\D{\mathscr{E}}.
}
$$
Let $\Theta:\D{\mathscr{E}}\to\D{\mathscr{E}}$ be the composition of all functors in the bottom row of the diagram. Then the restriction $\Theta^-$
of $\Theta$ to $\Df{\mathscr{E}}$ is the composition of the functors in the top row of the diagram. Since $\overline{D^-(i)}$ is an equivalence, $\Theta^-$
is naturally isomorphic to the identity functor of $\Df{\mathscr{E}}$ by $(\ast\ast)$.

Further, $\Theta$ commutes with direct products. This is obtained by showing that each component of $\Theta$ commutes with direct products.
In fact, since $\mathscr{E}$ is closed under products in $B\Modcat$ by Lemma \ref{thick}(1), $D(i)$ commutes with direct products. Clearly, the inclusion $\kappa$ preserves direct products. As $\Img(H)$ is a triangulated subcategory of $\D{B}$ closed under products by Theorem \ref{Torsion}, it follows from the dual result of \cite[Lemma 1.5]{BN} that $Q$ commutes with direct products. Finally, $\overline{\Phi}[-1]$ commutes with direct products by Lemma \ref{FF}(3).

Next, we show that $\Theta$ is dense, and therefore $\overline{\Phi}[-1]$ is  dense.

Let $\cpx{X}\in\C{\mathscr{E}}$. Then $X^{\leq n}\in\Cf{\mathscr{E}}$ for $n\geq 0$ and $\cpx{X}\simeq\underleftarrow{\holim}(X^{\leq n})$ in $\K{B}$ by Lemma \ref{Homotopy}. Recall that $\K{\mathscr{E}}$ is a triangulated subcategory of $\K{B}$ closed under direct products by Lemma \ref{thick}(1).
This implies $\cpx{X}\simeq\underleftarrow{\holim}(X^{\leq n})$ in $\K{\mathscr{E}}$ and thus also in $\D{\mathscr{E}}$ because the localization functor $\K{\mathscr{E}}\to\D{\mathscr{E}}$ preserves direct products.
As $\Theta$ commutes with products, it commutes with homotopy limits in $\D{\mathscr{E}}$.
Consequently, $\Theta(\cpx{X})=\Theta\big(\underleftarrow{\holim}(X^{\leq n})\big)\simeq\underleftarrow{\holim}\big(\Theta(X^{\leq n})\big)$.
Since the restriction of $\Theta$ to $\Df{\mathscr{E}}$ is isomorphic to the identity functor of $\Df{\mathscr{E}}$,
it follows from $X^{\leq n}\in\Df{\mathscr{E}}$ that $\underleftarrow{\holim}\big(\Theta(X^{\leq n})\big)\simeq \underleftarrow{\holim}(X^{\leq n})\simeq\cpx{X}$ in $\D{\mathscr{E}}$. Thus $\Theta(\cpx{X})\simeq\cpx{X}$ in $\D{\mathscr{E}}$. This shows that $\Theta$ is dense.

Finally, we show that $\overline{D(i)}$ is a triangle equivalence.

By Lemma \ref{1-tilting}(3), the functor $\overline{D(i)}$ is dense. By Lemma \ref{FF}(4), the composition of $\overline{\Phi}[-1]$ with $\overline{D(i)}$ is a triangle equivalence. Since $\overline{\Phi}[-1]$ is dense, $\overline{D(i)}$ is full. Moreover, by Lemma \ref{FF}(2), $\Ker(\overline{D(i)})=0$. In other words, $\overline{D(i)}$ sends nonzero objects of $\D{\mathscr{E}}$ to nonzero objects of $\Ker(G)$. It is known in \cite[p.446]{Ri} that a triangle functor between triangulated categories is faithful if it is full and sends nonzero objects to nonzero objects. This implies that $\overline{D(i)}$ is faithful, and thus an equivalence. $\square$

\subsection{Proofs of main result and its corollaries \label{sect3.3}}

We keep all notation introduced in the previous sections.

{\bf Proof of Theorem \ref{Main result}}. Since $\overline{D(i)}$ is a triangle equivalence by Proposition \ref{Embedding}, we can replace the left term $\Ker(G)$ of the recollement in Theorem \ref{Torsion} with $\D{\mathscr{E}}$ and obtain the recollement of the form
$$
\xymatrix@C=1.2cm{\D{\mathscr{E}}\ar[r]^-{D(i)}
&\D{B}\ar[r]^{G}\ar@/^1.2pc/[l]^{\Phi[-1]}\ar@/_1.2pc/[l]
&\D{A}\ar@/^1.2pc/[l]^H\ar@/_1.2pc/[l]_{j_!}}.$$
In this recollement, the functors $D(i)$, $\Phi[-1]$, $H$ and $G$ can be restricted to functors between $\mathscr{D}^{-}$-derived categories by Remark \ref{Bounded-above}. So, to show that this recollement restricts to $\mathscr{D}^{-}$-derived categories, it suffices to show that the left adjoint $j_!:\D{A}\to \D{B}$ of $G$ restricts to a functor $\Df{A}\to\Df{B}$.

For any $m\in\mathbb{Z}$ and $X\in \D{A}$, there are isomorphisms of abelian groups:
$$(\dag)\quad \Hom_{\D A}(X, G(B^\vee)[m])\simeq \Hom_{\D B}(j_!(X), B^\vee[m])\simeq
\Hom_{\K B}(j_!(X), B^\vee[m])$$
$$\;\;\quad\simeq H^m(j_!(X)^\vee)\simeq (H^{-m}(j_!(X)))^\vee.$$
Let $X\in\Df{A}$. Since $G(B^\vee)\in\Db{A}$ by Lemma \ref{Injective}, there is an integer $m_X$ such that $\Hom_{\D A}(X,$ $ G(B^\vee)[m])=0$ for all $m\leq m_X$.
Observe that the functor $(-)^{\vee}$ sends nonzero modules to nonzero modules. Thus $H^{-m}(j_!(X))=0$ for all $m\leq m_X$. This implies $j_!(X)\in\Df{B}$. $\square$

\medskip
Recall that a ring homomorphism $\lambda: B\to C$ is called a \emph{ring epimorphism} if the restriction functor $\lambda_*: C\Modcat\to B\Modcat$ is fully faithful.
A ring epimorphism $\lambda: B\to C$ is said to be \emph{homological} if $\Tor^B_m(C, C)=0$ for all $m\geq 1$. This is also equivalent to saying that the induced derived functor $D(\lambda_*): \D{C}\to \D{B}$ is fully faithful. Following \cite{XC2}, a good tilting module ${_A}T$ is  said to be \emph{homological} if there exists a homological ring epimorphism $\lambda:B\to C$ such that $D(\lambda_*)$ induces an equivalence $\D{C}\lraf{\simeq}\Ker(G)$ of triangulated categories.

To prove  Corollary \ref{Homological}, we recall the definition of derived decompositions of abelian categories in \cite{XC5}.

\begin{Def}{\rm \cite[Definition 1.1]{XC5}} \label{DC}
Let $\mathcal{A}$ be an abelian category, and let $\mathcal{B}$ and $\mathcal{C}$ be full subcategories of $\mathcal{A}$. The pair $(\mathcal{B}, \mathcal{C})$ is called a \emph{derived decomposition} of $\mathcal{A}$ if both $\mathcal{B}$ and $\mathcal{C}$ are abelian subcategories of $\mathcal{A}$ such that

$(1)$ the inclusions $\mathcal{B}\subseteq\mathcal{A}$ and $\mathcal{C}\subseteq\mathcal{A}$ induce fully faithful functors $\Db{\mathcal{B}}\hookrightarrow \Db{\mathcal{A}}$ and $\Db{\mathcal{C}}\hookrightarrow\Db{\mathcal{A}}$, respectively.

 $(2)$ $\big(\Db{\mathcal{B}},\Db{\mathcal{C}}\big)$ is a semi-orthogonal decomposition of $\Db{\mathcal{A}}$.
\end{Def}

{\bf Proof of Corollary \ref{Homological}}.  The equivalence of (1) and (2) is shown in \cite[Proposition 6.2]{BaPa},  while the one of (1) and (3) is proved in \cite[Theorem 1.1]{XC2}. Here, we prefer to giving a new and shorter proofs of these facts by the results in the previuos sections.

 $(1)\Rightarrow (2)$: If $T$ is homological, then, by definition, the associated homological ring epimorphism $\lambda:B\to C$
 induces an equivalence from $C\Modcat$ to $\mathscr{E}$. Thus $\mathscr{E}$ is an abelian subcategory of $B\Modcat$.

$(2)\Rightarrow (1)$: Suppose $(2)$ holds. Then $\mathscr{E}$ is closed under direct sums, direct products, kernels and cokernels by Lemma \ref{thick}(1).
It follows that there exists a ring epimorphism $\lambda: B\to C$  such that $\lambda_*$
induces an equivalence $C\Modcat\lraf{\simeq}\mathscr{E}$ of abelian categories (for example, see \cite[Lemma 2.1]{XC1}).
Since $\overline{D(i)}$ is fully faithful by Proposition \ref{Embedding}, the functor $D(\lambda_*)$ is also fully faithful.
Thus $\lambda$ is homological. Further, it follows from $\Img(\overline{D(i)})=\Ker(G)$ that $D(\lambda_*)$ induces an equivalence from $\D{C}$ to $\Ker(G)$.
This shows $(1)$.

$(2)\Rightarrow (3)$: Let $\cpx{\rho}: \Hom_A(\cpx{P}, A)\otimes_AT\to \Hom_A(\cpx{P}, T)$ be the chain map defined in Lemma \ref{Injective}(2).
Associated with $B^\vee$, there exists the triangle in $\D{B}:$
$$(\bigtriangleup): \; \cpx{X}[-1]\to B^\vee\to \cpx{Y}\to \cpx{X}$$ with
$\cpx{Y}:=(\Hom_A(\cpx{P}, A)\otimes_AT)^\vee\in\Img(H)=\Img(D(j))\;\;\mbox{and}\;\; \cpx{X}:=\Coker(\cpx{\rho})^\vee[1]\in\Ker(G)=\Img(D(i)).$
Since $\mathcal{Y}\subseteq B\Modcat$ is an abelian subcategory, we have $H^s(\cpx{Y})\in\mathcal{Y}$ for any $s\in\mathbb{Z}$. Taking cohomologies on the triangle yields $H^{-s}(\cpx{X})\simeq H^{-s}(\cpx{Y})$ whenever $s\geq 2$.
Moreover, $H^t(\cpx{X})\simeq H^{t+1}\big(\Coker(\cpx{\rho})^\vee\big)\simeq \big(H^{-t-1}(\Coker(\cpx{\rho}))\big)^\vee=0$ for any $t\geq 0$. In particular, $H^0(\cpx{X})=0$.

Suppose $(2)$ holds. Then $H^s(\cpx{X})\in\mathscr{E}$ for any $s\in\mathbb{Z}$.
By Lemma \ref{Ext-orthogonal}, $H^{-s}(\cpx{X})=0=H^{-s}(\cpx{Y})$ for any $s\geq 2$. It follows that
$0=H^{-m}(\cpx{Y})\simeq (H^m(\Hom_A(\cpx{P}, A)\otimes_AT))^\vee$ for any $m\geq 2$.
Now $(3)$ holds by the fact that the functor $(-)^{\vee}$ sends nonzero modules to nonzero modules.

$(3)\Rightarrow (2)$:  Suppose $(3)$ holds. By Lemma \ref{thick}(1), $\mathscr{E}$ is an abelian subcategory of $B\Modcat$ if and only if
it is closed under cokernels in $B\Modcat$. Recall that $\mathcal{X}:=\{X\in B\Modcat\mid \Ext_B^n(X,Y)=0,\, Y\in\mathcal{Y},\, n=0,1 \}$ is always closed under cokernels in $B\Modcat$. So it is enough to show
$\mathscr{E}=\mathcal{X}$.

Clearly, $\mathscr{E}\subseteq \mathcal{X}$. To show the converse, we first claim that each module $U\in\mathcal{X}$ can be embedded into a module belonging to $\mathscr{E}$. Actually, $U$ can be embedded into its injective envelop $J$ in $B\Modcat$. By the first isomorphism in Corollary \ref{adjoint}(2), $U$ is isomorphic to a submodule of $X_J$.  So, it suffices to show  $X_I\in\mathscr{E}$ for any injective $B$-module $I$.

Obviously, $H^s(\cpx{Y})=0$ for any $s\geq 1$. Combining this with $(3)$, we have $H^{-s}(\cpx{Y})=0$ for any $s\neq 0, 1$. This implies
$H^{-t}(\cpx{X})=0$ for any $t\neq 0, 1$. As $H^0(\cpx{X})=0$, we get $\cpx{X}\simeq H^{-1}(\cpx{X})[1]$ in $\D{B}$.
Consequently, $H^{-1}(\cpx{X})\in\mathscr{E}$ by $\cpx{X}\in\Ker(G)$. By Lemmas \ref{1-tilting}(2) and \ref{Injective}(2), the triangle $(\diamondsuit)$ associated with $B^\vee$ is isomorphic to the one $(\bigtriangleup)$. This shows $X_{B^\vee}\simeq H^{-1}(\cpx{X})\in\mathscr{E}$. Recall that $D(i)$ and $\Phi$ commute with direct products by Lemma \ref{FF}(3), and that $\Prod(B^\vee)$ consists of all injective $B$-modules.
Thus $X_I\in\Prod(X_{B^\vee})\subseteq \mathscr{E}$ by the diagram $(\sharp)$ and Lemma \ref{thick}(1).

Now, it follows from the above claim that there is an injection $f_0:U\to E_0$ with $E_0\in\mathscr{E}$. Since $\mathcal{X}$ contains $\mathscr{E}$ and is closed under cokernels in $B\Modcat$, we have $\Coker(f_0)\in\mathcal{X}$. By iterating this construction, we can get an infinitely long exact sequence
$0\to U\lraf{f_0} E_0\to E_1\to \cdots\to E_i\to\cdots $ of $B$-modules such that $E_i\in\mathscr{E}$ for all $i\in\mathbb{N}$.
Then $U\in\mathscr{E}$ by Lemma \ref{thick}(2). Hence we have shown $\mathcal{X}\subseteq \mathscr{E}$. Thus $(2)$ holds.

(2) $\Leftrightarrow$ (4).
Note that $D^b(i):\Db{\mathscr{E}}\to\Db{B}$ and $D^b(j): \Db{\mathcal{Y}}\to \Db{B}$ are fully faithful by Proposition \ref{Embedding} and Lemma \ref{FF}(1), respectively. Since ${_A}T$ and $T_B$ have finite projective dimension, $H$ and $G$ can be restricted to functors between bounded derived categories. It follows that $(\Db{B}\cap\Ker(G), \Db{B}\cap\Img(H))$ is a semi-orthogonal decomposition of $\Db{B}$. Further, $\Img(\Db{i})=\Db{B}\cap\Ker(G)$ by Lemmas \ref{1-tilting}(3) and \ref{thick}(2), while $\Img(\Db{j})=\Db{B}\cap\Img(H)$ by Lemma \ref{1-tilting}(3) and Corollary \ref{Restriction}(2).
Thus $\big(\Img(\Db{i}), \Img(\Db{j})\big)$ is a semi-orthogonal decomposition of $\Db{B}$. This implies that $(\mathscr{E}, \mathcal{Y})$ is a derived decomposition of $B\Modcat$ if and only if $\mathscr{E}$ is an abelian subcategory of $B\Modcat$.
In other words, $(2)$ and $(4)$ are equivalent.
$\square$

\medskip
Combining Corollary \ref{Homological} with \cite[Corollary 1.3]{XC5}, we have the following result where our new contribution is
about the bound of projective dimension.
\begin{Koro}
If ${_A}T$ is a homological tilting module, then there exists a homological ring epimorphism $\lambda: B\to C$  such that
$\pd({_B}C)\leq 1$ and $\lambda_*$ induces an equivalence from $C\Modcat$ to $\mathscr{E}$.
\end{Koro}

{\bf Proof of Corollary \ref{Koro1.4}}. By Lemma \ref{thick}(2), up to triangle equivalence, we have the identifications of categories:
$\mathscr{D}^b(\mathscr{E})=\{\cpx{X}\in \D{\mathscr{E}}\mid H^i(\cpx{X})=0, |i| >>0 \}$, $\mathscr{D}^+(\mathscr{E})=\{\cpx{X}\in \D{\mathscr{E}}\mid H^i(\cpx{X})=0, i<<0 \}$ and $\mathscr{D}^-(\mathscr{E})=\{\cpx{X}\in \D{\mathscr{E}}\mid H^i(\cpx{X})=0, i>>0 \}$. Since $\pd(_AT)\le n$ and $\pd(T_B)\le n$,
the functors $G$ and $H$ in Theorem \ref{Main result} restrict to functors at $\mathscr{D}^*$-level for $*\in \{b,+,-\}$. So we get the half recollement
$$\xymatrix{\mathscr{D}^*(\mathscr{E})\ar[r]&\mathscr{D}^*(B)\ar[r]^{G}\ar@/^1.2pc/[l]
&\mathscr{D}^*(A)\ar@/^1.2pc/[l]_{H}}\vspace{0.2cm}$$
To show Corollary \ref{Koro1.4}, it suffices to show that the left adjoint $j_!:\D{A}\to \D{B}$ of $G:\D{B}\to\D{A}$ restricts to a functor $\mathscr{D}^*(A)\to\mathscr{D}^*(B)$ whenever $A$ is left coherent.

Since $A$ is left coherent, the direct products of projective $A\opp$-modules are flat. As $(-)^{\vee}:A\opp\Modcat\to A\Modcat$ sends flat modules to injective modules, $\Hom_A(\cpx{P}, A)^\vee$ is a bounded complex of injective $A$-modules. By Lemma \ref{Injective}(1) and the proof of Theorem \ref{Main result} (see $(\dag))$, there are
isomorphisms {\small $$(H^{-m}(j_!(X)))^\vee\simeq \Hom_{\D A}(X, G(B^\vee)[m])\simeq
\Hom_{\D A}(X, \Hom_A(\cpx{P}, A)^\vee[m])\simeq\Hom_{\K A}(X, \Hom_A(\cpx{P}, A)^\vee[m]) $$}for any $X\in\D{A}$ and $m\in\mathbb{Z}$.
Note that the functor $(-)^{\vee}$ sends nonzero modules to nonzero modules. Thus $X\in \mathscr{D}^*(A)$ if and only if $j_!(X)\in\mathscr{D}^*(B)$. Further, the right (left) adjoint functor from $\D{B}$ to $\D{\mathscr{E}}$ of $D(i)$ restricts to $\mathscr{D}^*$-level follows from the first (second) triangle in Definition \ref{def01}(4). $\square$

{\footnotesize
\bigskip Hongxing Chen, School of Mathematical Sciences  \&  Academy for Multidisciplinary Studies, Capital Normal University, 100048 Beijing, P. R. China

{\tt Email: chenhx@cnu.edu.cn}

\bigskip
Changchang Xi, School of Mathematical Sciences, Capital Normal
University, 100048 Beijing, China, and

School of Mathematics and Information Science,
Henan Normal University, Xinxiang, Henan, China

{\tt Email: xicc@cnu.edu.cn} }

{\url{http://math0.bnu.edu.cn/~ccxi/}}


{\footnotesize
\begin{thebibliography}{99}
\bibitem{AC} {{\sc L. Angeleri-H\"{u}gel} and  {\sc F. U. Coelho,}
Infinitely generated tilting modules of finite projective dimension.
{\it Forum Math.} \textbf{13} (2001) 239-250.}

\bibitem{HKL}{{\sc L. Angeleri-H\"ugel, S. K\"onig} and {\sc Q. H. Liu}, Recollements
and tilting objects, {\it J. Pure Appl. Algebra} \textbf{215} (2011)
420-438}.

\bibitem{Bz1}{{\sc S. Bazzoni},
The t-structure induced by an $n$-tilting module. \emph{Tran.
Amer. Math. Soc.} \textbf{371} (2019) 6309-6340.}

\bibitem{Bz2}{{\sc S. Bazzoni, F. Mantese} and {\sc A. Tonolo},
Derived equivalence induced by $n$-tilting modules. \emph{Proc.
Amer. Math. Soc.} \textbf{139} (2011) 4225-4234.}

\bibitem{BaPa}{{\sc S. Bazzoni} and {\sc A. Pavarin},
Recollements from partial tilting complexes. \emph{J. Algebra}
\textbf{388} (2013) 338-363.}

\bibitem{BBD}{{\sc A. A. Beilinson, J. Bernstein} and {\sc P.
Deligne},  Faisceaux pervers, \emph{Asterisque} \textbf{100} (1982)
5-171.}

\bibitem{BI}{{\sc A. Beligiannis} and {\sc I. Reiten},
Homological and homotopical aspects of torsion theories, {\it Mem.
Amer. Math. Soc}. {\bf 188} (2007), no. 883, 1-207.}

\bibitem{BN}{{\sc M. B\"okstedt} and {\sc A. Neeman},
Homotopy limits in triangulated categories, \emph{Compositio Math.}
\textbf{86} (1993) 209-234.}

\bibitem{BBu}{{\sc S. Brenner} and {\sc M. R. Butler}, Generalizations of the Bernstein-Gelfand-Ponomarev reflection functors’,
Representation theory II, Springer Lecture Notes in Mathematics \textbf{83}2 (eds V. Dlab and P. Gabriel; Springer,
Berlin, 1980) 103-169.}

\bibitem{XC1}{{\sc H. X. Chen} and {\sc C. C. Xi}, Good tilting
modules and recollements of derived module categories, \emph{Proc.
Lond. Math. Soc.} \textbf{104} (2012) 959-996.}

\bibitem{XC2}{{\sc H. X. Chen} and {\sc C. C. Xi}, Good tilting
modules and recollements of derived module categories, II.
\emph{J. Math. Soc. Japan} \textbf{71}(2) (2019) 515-554.}

\bibitem{XC4}{{\sc H. X. Chen} and {\sc C. C. Xi}, Recollements of derived categories,
III: Finitistic dimensions. J. London Math. Soc. (2) \textbf{95} (2017) 633-658.}


\bibitem{XC5}{{\sc H. X. Chen} and {\sc C. C. Xi}, Derived decompositions of abelian categories, I.
\emph{Pacific J. Math.} (to appear). Primary version at arXiv:1804.10759.}

\bibitem{CPS}{{\sc E. Cline, B. Parshall} and {\sc L. Scott},
Derived categories and Morita theory, \emph{J. Algebra} \textbf{104}
(1986) 397-409.}

\bibitem{CPS2}{{\sc E. Cline, B. Parshall} and {\sc L. Scott},
Algebraic stratification in representation categories, \emph{J.
Algebra} \textbf{117} (1988) 504-521.}

\bibitem{ct}{{\sc R. Colpi} and {\sc J. Trlifaj}, Tilting modules and tilting torsion
theories, \emph{J. Algebra} \textbf{178} (1995) 614-634.}

\bibitem{happel87}{{\sc D. Happel}, On the derived category of a finite-dimensional algebra,
\emph{Comment. Math. Helv.} \textbf{62}(3) (1987) 339-389.}

\bibitem{hr}{{\sc D. Happel} and {\sc C. M. Ringel}, Tilted algebras.
\emph{Trans. Amer. Math. Soc.} \textbf{274} (1982) 399-443.}

\bibitem{happel2}{{\sc D. Happel}, Reduction techniques for homological conjectures. \emph{Tsukuba
J. Math.} \textbf{17} (1993) 115-130.}

\bibitem{Huybrechts}{ {\sc  D. Huybrechts}, \emph{Fourier-Mukai transforms in algebraic geometry},
Oxford Math. Monogr., The Clarendon Press/Oxford University Press, Oxford, 2006.}

\bibitem{Keller}{{\sc B. Keller}, Derived categories and their uses. In: M. Hazewinkel, Ed., Handbook of algebra, vol. 1, (1996) 671-701.}

\bibitem{mntty}{{\sc H. Matsui, T. T. Nam, R. Takahashi, N. M. Tri} and {\sc D. N. Yen}, Cohomological dimensions of specialization-closed subsets and subcategories of modules. \emph{Proc. Amer. Math. Soc.}  \textbf{149}(2) (2021) 481-496.}

\bibitem{miy}{{\sc Y. Miyashita}, Tilting modules of finite projective
dimension, \emph{Math. Z.} \textbf{193} (1986) 113-146.}


\bibitem{Neeman} {{\sc A. Neeman}, {\it  Triangulated categories},
Ann. of Math. Stud. {\bf 148}, Princeton University Press, 2001.}

\bibitem{Neeman1}{{\sc  A. Neeman}, The homotopy category of flat modules, and Grothendieck duality, \emph{Invent. Math.}
\textbf{174} (2008) 255-308.}

\bibitem{Pos1}{{\sc L. Positselski}, Contraherent cosheaves, arXiv:1209.2995v6 (2017).}

\bibitem{Ri}{{\sc J. Rickard}, Morita theory for
derived categories. {\it J. London Math. Soc.}
\textbf{39} (1989) 436-456.}

\bibitem{stenstroem}{{\sc B. Stenstr\"om}, \emph{Rings of quotients}: An introduction to methods of ring theory. Springer-Verlag
Berlin Heidelberg New York, 1975.}

\bibitem{STK}{{\sc J. \v{S}\v{t}ov\'{i}\v{c}ek}, Derived equivalences induced by big cotilting modules,
{\it Adv. Math.} \textbf{263} (2014) 45-87.}

\bibitem{We}{{\sc C. A. Weibel}, {\it An introduction to homological
algebra}, Camb. Studies in Adv. Math. {\bf 38}, Cambridge University
Press, 1994.}

\bibitem{yang} {{\sc D. Yang}, Recollements from generalized tilting, \emph{Proc. Amer. Math. Soc.} {\bf 140} (2012) 83-91.}

\end{thebibliography}}
\end{document}